\documentclass[amscd,amssymb]{amsart}

\usepackage{amssymb}
\usepackage[all]{xy}

\newtheorem{thm}{Theorem}[section]
\newtheorem{lem}[thm]{Lemma}
\newtheorem{cor}[thm]{Corollary}
\newtheorem{prop}[thm]{Proposition}
\newtheorem{conj}[thm]{Conjecture}

\theoremstyle{definition}

\newtheorem{defn}[thm]{Definition}

\newtheorem{ex}[thm]{Example}
\newtheorem{exs}[thm]{Examples}

\theoremstyle{remark}

\newtheorem{rem}[thm]{Remark}

\numberwithin{equation}{section}

\newcommand{\thmref}[1]{Theorem~\ref{#1}}
\newcommand{\corref}[1]{Corollary~\ref{#1}}
\newcommand{\secref}[1]{\S\ref{#1}}
\newcommand{\conjref}[1]{Conjecture~\ref{#1}}
\newcommand{\propref}[1]{Proposition~\ref{#1}}
\newcommand{\lemref}[1]{Lemma~\ref{#1}}

\newcommand{\exref}[1]{Example~\ref{#1}}

\newcommand{\colim}{\operatorname*{colim}}

\newcommand{\Hom}{\operatorname{Hom}}
\newcommand{\Ext}{\operatorname{Ext}}

\newcommand{\Ind}{\operatorname{Ind}}
\newcommand{\Res}{\operatorname{Res}}

\newcommand{\coker}{\operatorname{coker}}
\newcommand{\im}{\operatorname{im}}
\newcommand{\Inf}{\operatorname{Inf}}
\newcommand{\Rep}{\operatorname{Rep}}

\newcommand{\A}{{\mathcal  A}}

\newcommand{\U}{{\mathcal  U}}
\newcommand{\K}{{\mathcal  K}}

\newcommand{\PP}{{\mathcal P}}

\newcommand{\Z}{{\mathbb  Z}}

\newcommand{\F}{{\mathbb  F}}

\newcommand{\ra}{\rightarrow}
\newcommand{\xra}{\xrightarrow}

\newcommand{\xla}{\xleftarrow}

\newcommand{\hra}{\hookrightarrow}
\newcommand{\era}{\twoheadrightarrow}
\newcommand{\xera}[1]{\overset{#1}{\twoheadrightarrow}}

\begin{document}

\title[group cohomology primitives]{Primitives and central detection numbers in group cohomology}

 \author[Kuhn]{Nicholas J.~Kuhn}
 \address{Department of Mathematics \\ University of Virginia \\ Charlottesville, VA 22904}
 \email{njk4x@virginia.edu}
\thanks{This research was partially supported by a grant from the National Science Foundation}

 \date{December 5, 2006.}


 \begin{abstract}
Fix a prime $p$. Given a finite group $G$, let $H^*(G)$ denote its mod $p$ cohomology.
In the early 1990's, Henn, Lannes, and Schwartz introduced two invariants $d_0(G)$ and $d_1(G)$ of $H^*(G)$ viewed as a module over the mod $p$--Steenrod algebra.  They showed that, in a precise sense, $H^*(G)$ is respectively detected and determined by $H^d(C_G(V))$ for $d \leq d_0(G)$ and $d \leq d_1(G)$, with $V$ running through the elementary abelian $p$--subgroups of $G$.

The main goal of this paper is to study how to calculate these invariants.  We find that a critical role is played by the image of the restriction of $H^*(G)$ to $H^*(C)$, where $C$ is the maximal central elementary abelian $p$--subgroup of $G$.  A measure of this is the top degree $e(G)$ of the finite dimensional Hopf algebra $H^*(C)\otimes_{H^*(G)}\F_p$, a number that tends to be quite easy to calculate.

Our results are complete when $G$ has a $p$--Sylow subgroup $P$ in which every element of order $p$ is central.  Using Benson--Carlson duality, we show that in this case, $d_0(G) = d_0(P) = e(P)$, and a similar exact formula holds for $d_1$.  As a bonus, we learn that $H^{e(G)}(P)$ contains nontrivial essential cohomology, reproving and sharpening a theorem of Adem and Karagueuzian.

In general, we are able to show that $d_0(G) \leq \max \{e(C_G(V)) \ | \ V < G\}$ if certain cases of Benson's Regularity Conjecture hold.  In particular, this inequality holds for all groups such that the difference between the $p$--rank of $G$ and the depth of $H^*(G)$ is at most 2.  When we look at examples with $p=2$, we learn that $d_0(G) \leq 14$ for all groups with 2--Sylow subgroup of order up to 64, with equality realized when $G = SU(3,4)$.

Enroute we study two objects of independent interest.
If $C$ is any central elementary abelian $p$--subgroup of $G$, then $H^*(G)$ is a $H^*(C)$--comodule, and we prove that the subalgebra of $H^*(C)$--primitives is always Noetherian of Krull dimension equal to the $p$--rank of $G$ minus the $p$--rank of $C$.
If the depth of $H^*(G)$ equals the rank of $Z(G)$, we show that the depth essential cohomology of $G$ is nonzero (reproving and extending a theorem of Green), and Cohen--Macauley in a certain sense, and prove related structural results.

\end{abstract}

\maketitle

\section{Introduction}

Fix a prime $p$, and let $H^*(G)$ denote the mod $p$ cohomology ring of a finite group $G$.  The $p$--elementary abelian subgroups of $G$ have had a featured role in the study of group cohomology since D.Quillen's famous work \cite{quillen} in the late 1960's.  In particular, these subgroups become the objects in a category $\A(G)$ having morphisms the homomorphisms generated by subgroup inclusion and conjugation by elements in $G$.  The inclusions $V < G$ then induce a map
$$  H^*(G) \xra{\lambda_0} \lim_{V \in \A(G)} H^*(V) \subseteq \prod_V H^*(V), $$
and $\lambda_0$ is shown to have kernel and cokernel that are nilpotent in an appropriate sense.

Viewing $H^*(G)$ as the mod $p$ cohomology of the classifying space $BG$ makes it evident that $H^*(G)$ is an object in $\K$ and $\U$, the categories of unstable algebras and modules over the mod $p$ Steenrod algebra $\A$.  The 1980's and 1990's saw a revolution in our understanding of these categories, and the 1995 paper  of H.-W.Henn, J.Lannes, and L.Schwartz \cite{hls1} revisited Quillen's approximation of $H^*(G)$ from this new perspective.

For each $d \geq 0$, the group homomorphisms $V \times C_G(V) \ra G$ induce a map of unstable algebras  $$\displaystyle H^*(G) \ra  \prod_V H^*(V) \otimes H^{\leq d}(C_G(V)),$$ where $M^{\leq d}$ denotes the quotient of a graded module $M^*$ by all elements of degree more than $d$.  The image of this map lands in an evident subalgebra of `compatible' elements which Henn, Lannes, and Schwartz show can be naturally identified with $L_dH^*(G)$, where $L_d: \U \ra \U$ is localization away from the localizing subcategory generated by $(d+1)$--fold suspensions of unstable modules.  Thus Quillen's map can be viewed as the just the bottom of a tower of localizations of $H^*(G)$ associated to the nilpotent filtration of $\U$:
\begin{equation*}
\xymatrix{
&&& \vdots \ar[d] \\
&&& L_2H^*(G) \ar[d]^{p_2} \\
&&& L_1H^*(G) \ar[d]^{p_1} \\
H^*(G) \ar[rrr]^{\lambda_0}  \ar[urrr]^{\lambda_1} \ar[uurrr]^{\lambda_2} &&& L_0H^*(G),
}
\end{equation*}
where we have
$$ H^*(G) \xra{\lambda_d} L_dH^*(G) \subseteq \prod_V H^*(V) \otimes H^{\leq d}(C_G(V)).$$

This caused the authors of \cite{hls1} to introduce two new invariants of $G$: $d_0(G)$ and $d_1(G)$ are the smallest $d$'s such that $H^*(G)$ is respectively detected by, and isomorphic to, $L_dH^*(G)$.  Alternatively, $d_0(G)$ is the smallest $d$ such that $H^*(G)$ contains no $(d+1)$--fold suspensions of a nontrivial unstable module, and $d_1(G)$ is the smallest $d$ such that also $\Ext_{\A}^1(\Sigma^{d+1}N,H^*(G))=0$ for all $N \in \U$.

These invariants satisfy a few easily verified nice properties: $d_0(G \times H) = d_0(G) + d_0(H)$, $d_1(G \times H) = \max \{ d_1(G)+ d_0(H),d_0(G) +d_1(H) \}$, and $d_i(G) \leq d_i(P)$ if $P$ is a $p$--Sylow subgroup of $G$.  However, they are not well behaved under taking subgroups, quotient groups, and extensions; e.g., every $G$ embeds in a symmetric group $\Sigma_n$ and $d_0(\Sigma_n) = 0$.  Rough upper bounds for $d_0(G)$ and $d_1(G)$ were found in \cite{hls1}; e.g., $d_0(G)$ is bounded by $n^2$ if a $p$--Sylow subgroup of $G$ admits a faithful $n$ dimensional complex representation. However, in all but a few examples, these bounds seem far from optimal.  Up to now, what determines these group invariants has remained mysterious, and they have not been connected to other work in group cohomology.

A main goal of this paper is to present a way to calculate the number $d_0(G)$, and, in some cases, $d_1(G)$.  Our finding is that these numbers seem to be controlled by the restriction of cohomology to maximal central $p$--elementary abelian subgroups.  Our results are complete when $G$ has a Sylow subgroup that is $p$--central, i.e. a group in which every element of order $p$ is central.  For example, when $p=2$, we compute that $d_0(SU(3,4)) = 14$ and $d_1(SU(3,4)) = 18$, where, by constrast, the estimates from \cite{hls1} yield only that $d_0(SU(3,4)) \leq 64$ and $d_1(SU(3,4)) \leq 120$.

Our method is to combine $\U$--technology, in the spirit of \cite{hls1}, with duality results as in the work of D.Benson and J.Carlson \cite{bc}.  We ultimately connect a conjectured upper bound for $d_0(G)$ to Benson's Regularity Conjecture \cite{b}, known to hold if the $p$--rank of $G$ and the $p$--rank of $Z(G)$ differ by at most 2.  This is the case for all 2-groups of order 64 or less, and, using cohomology calculations from \cite{carlson et al}, we've been able to verify by hand that $d_0(G) \leq 14$ for all such groups.

A number of side results of independent interest come up in our investigations.

We are led to study carefully the cohomology of central extensions, in particular the structure of associated algebras of primitives.  One outcome of this is a new proof of A.Adem and D.Karagueuzian's theorem \cite{ak} that $p$--central $p$--groups have nonzero essential cohomology.  We show that in an explicit degree there is a nonzero cohomology class that is simultaneously essential and annihilated by all Steenrod operations of positive degree.

Deriving our general estimate of $d_0(G)$ involves a careful study of the depth essential cohomology of Carlson, et. al. \cite{carlson et al} in the important special case that the depth of $H^*(G)$ equals the rank of the center.  We prove that then the depth essential cohomology is both nonzero - reproving the main theorem of \cite{green1} without D.Green's hypothesis that $G$ be a $p$--group - and Cohen--Macauley.

In the next section we describe our results in more detail.

\section{Main Results} \label{main results section}

\subsection{The cohomology of central extensions}

Suppose we have a central extension of finite groups
$$ C \xra{i} G \xra{q} Q,$$
where $C$ is $p$--elementary abelian of rank $c$.

We define various objects associated to this situation.

The extension corresponds to an element $\tau \in H^2(Q;C)$. Since $H^2(Q;C) = \Hom(H_2(Q),C)$, the extension can also be considered as corresponding to a homomorphism $\tau:H_2(Q) \ra C$, or, equivalently, its dual $\tau^{\#}: C^{\#} \ra H^2(Q)$.

Let $\{E_r^{*,*}\}$ denote the Serre spectral sequence associated to the extension, converging to $H^*(G)$, and with $E_2^{*.*} = H^*(Q) \otimes H^*(C)$.  Under the identification $C^{\#} = H^1(C)$, it is standard that $\tau^{\#}$ corresponds to $ d_2: E_2^{0,1} \ra E_2^{2,0}$.

Let $I_{\tau} \subset H^*(Q)$ be the ideal generated by $\A \cdot \im(\tau_G^{\#})$, so that $H^*(Q)/I_{\tau}$ is an unstable algebra.  It is easy to see that $I_{\tau}$ is contained in the kernel of inflation $q^*:H^*(Q) \ra H^*(G)$.

Call a subalgebra $A$ of $H^*(G)$ a {\em $(G,C)$--Duflot subalgebra}, if the composite $A \subseteq H^*(G) \xra{i^*} \im(i^*)$ is an isomorphism, where $i^*: H^*(G) \ra H^*(C)$ is the restriction. As we will describe more precisely in \secref{res subsection}, as an algebra, the Hopf algebra $\im(i^*) \subseteq H^*(C)$ will necessarily be free graded commutative on $c$ polynomial generators, possibly tensored with an exterior algebra on some generators in degree 1, if $p$ is odd.  It follows that Duflot subalgebras exist and have the same form.  Let $Q_AH^*(G)$ denote the graded algebra\footnote{$Q_AH^*(G)$ will not necessarily be an unstable algebra, as $A$ need not be closed under Steenrod operations.} of $A$--indecomposables $H^*(G) \otimes_A \F_p$, or, equivalently, the quotient of $H^*(G)$ by the ideal generated by the positive degree elements of $A$.

$H^*(C)$ is a Hopf algebra, and the multiplication map $m: C \times G \ra G$ induces a map of unstable algebras
$$m^*: H^*(G) \ra H^*(C) \otimes H^*(G)$$ making $H^*(G)$ into a $H^*(C)$--comodule.
We define the associated algebra of primitives to be
\begin{equation*}
\begin{split}
P_CH^*(G) & = \{ x \in H^*(G) \ | \ m^*(x) = 1 \otimes x \} \\
& = \text{Eq } \{ H^*(G)
\begin{array}{c} m^* \\[-.08in] \longrightarrow \\[-.1in] \longrightarrow \\[-.1in] \pi^*
\end{array}
H^*(C \times G) \},
\end{split}
\end{equation*}
where $\pi: C \times G \ra G$ is the projection.
It is easy to check that $P_CH^*(G)$ is an unstable algebra that contains the image of the inflation map.  Thus $q^*: H^*(Q) \ra H^*(G)$ refines to a map of unstable algebras
$$ q_{\tau}: H^*(Q)/I_{\tau} \ra P_CH^*(G).$$

\begin{thm} \label{prim thm} With the notation as above, the following are true.  \\

\noindent{\bf (a)} $H^*(G)$ is a free $A$--module. Moreover $\{E_r^{*,*}\}$ is a spectral sequence of free $E_{\infty}^{0,*}$--modules, and applying $Q_{E_{\infty}^{0,*}}$ to the spectral sequence yields a spectral sequence converging to $Q_AH^*(G)$ with $E_2$--term $Q_{E_{\infty}^{0,*}}H^*(C) \otimes H^*(Q)$.\\

\noindent{\bf (b)} The composite $P_CH^*(G) \hra H^*(G) \era Q_AH^*(G)$ is monic. \\

\noindent{\bf (c)}  Both $P_CH^*(G)$ and $Q_AH^*(G)$ are finitely generated $H^*(Q)$--modules. \\

\noindent{\bf (d)}  The map $q_{\tau}: H^*(Q)/I_{\tau} \ra P_CH^*(G)$ is an $F$--isomorphism\footnote{In the sense of Quillen \cite{quillen}: $\ker(\bar{q}^*)$ is nilpotent, and for all $x \in P_CH^*(G)$, there exists a $k$ so that $x^{p^k} \in \im(\bar{q}^*)$.}, and the rings $H^*(Q)/I_{\tau}$,  $\im(q^*)$, $P_CH^*(G)$, and $Q_AH^*(G)$, are all Noetherian of Krull dimension equal to (the $p$--rank of $G$) - (the rank of $C$).
\end{thm}

Let $C(G) <G$ be the $p$--elementary abelian part of $Z(G)$.  If $C=C(G)$, the first part of statement (a) recovers J. Duflot's result \cite{duflot} that the depth of $H^*(G)$ is at least as great as the rank of $C(G)$.\footnote{We are claiming no originality in the proof of this, which is similar to all proofs of Duflot's theorem following \cite{broto henn}.  The spectral sequence refinement seems to be a new observation.} We will call a $(G,C(G))$--Duflot subalgebra of $H^*(G)$ simply a {\em Duflot subalgebra}.

The $p$--rank of $G$ equals the rank of $C$ exactly when $C=C(G)$ and $G$ is $p$--central, and we have the following corollary.

\begin{cor} \label{PC prim cor}  If $G$ is $p$--central, and $C=C(G)$, then the rings $H^*(Q)/I_{\tau}$, $\im(q^*)$, $P_{C}H^*(G)$, and $Q_AH^*(G)$ all have Krull dimension zero and so are finite dimensional $\F_p$--algebras.
\end{cor}

\subsection{Quillen's category and functors involving primitives}  Our algebras of primitives arise in two formulae associated to $H^*(G)$, viewed as an object in $\K$.  To describe these, we need to introduce some notation.

Given a small category $\mathcal C$, we let $\mathcal C^{\#}$ denote the associated twisted arrow category:
the objects of $\mathcal C^{\#}$ are the morphisms of $\mathcal C$, and a morphism $\alpha \rightsquigarrow \beta$ from $\alpha: A_1 \ra A_2$ to $\beta: B_1 \ra B_2$ is a commutative diagram in $\mathcal C$
\begin{equation*}
\xymatrix{
A_1  \ar[d] \ar[r]^{\alpha} & A_2   \\
B_1 \ar[r]^{\beta} & B_2. \ar[u] }
\end{equation*}

The functor assigning $H^*(V)$ to $V \in \A(G)$ is contravariant, while the assignment of $H^*(C_G(V))$ is covariant.  Now observe that the assignment of $P_{\alpha(V_1)}H^*(C_G(V_2))$ to $\alpha: V_1 \ra V_2$ can be viewed as defining a contravariant functor of $\A(G)^{\#}$.

Let $\A_C(G)$ denote the full subcategory of $\A(G)$ having as objects the $V$ containing $C(G)$.  If $G$ is $p$--central, then $\A_C(G)$ has a single object and morphism.

\subsection{A formula for the locally finite part of $H^*(G)$}

If $M$ is an unstable $\A$--module, we define $M_{LF}$, the locally finite part of $M$, by
$$ M_{LF} = \{ x \in M \ | \A x \subset M \text{ is finite}  \}.$$
This is again an unstable module, and is an unstable algebra if $M$ is.

\begin{thm} \label{LF thm} There is a natural isomorphism of unstable algebras $$\displaystyle H^*(G)_{LF} \simeq \lim_{V_1 \xra{\alpha} V_2} P_{\alpha(V_1)}H^*(C_{G}(V_2)),$$
where the limit is over $\A_C(G)^{\#}$.
\end{thm}

\begin{cor} \label{LF cor} If $G$ is $p$--central, then
$\displaystyle H^*(G)_{LF} =  P_{C(G)}H^*(G).$
\end{cor}

\subsection{A formula for $\bar R_dH^*(G)$}

An unstable module $M \in \U$ has a canonical `nilpotent' filtration \cite{s1,k2,hls1}:
$$ \dots \subseteq nil_2M \subseteq nil_1 M \subseteq nil_0M = M.$$
In general, $nil_dM/nil_{d+1}M = \Sigma^dR_dM$, where $R_dM$ is reduced, i.e. has no nontrivial submodules that are suspensions.  We let $\bar R_dM$ denote the nilclosure $L_0R_dM$ of $R_dM$.

The module $nil_dM$ identifies with the kernel of $\lambda_d: M \ra L_{d-1}M$, and a bit of diagram chasing will show that $\Sigma^d \bar R_dM$ is isomorphic to the kernel of $L_dM \ra L_{d-1}M$: see \propref{loc prop}. Thus $d_0(G)$ is the length of the filtration of $H^*(G)$, and also is the biggest $d$ such that $\bar R_dH^*(G) \neq 0$.

\begin{thm} \label{Rd thm} There is a natural isomorphism of unstable modules
$$\displaystyle \bar R_dH^*(G) \simeq \lim_{V_1 \xra{\alpha} V_2} H^*(V_1) \otimes P_{\alpha(V_1)}H^d(C_{G}(V_2)),$$
where the limit is over $\A_C(G)^{\#}$.
\end{thm}

\begin{cor} \label{Rd cor} If $G$ is $p$--central, then there is an isomorphism of unstable modules
$$\displaystyle \bar R_dH^*(G) \simeq  H^*(C(G)) \otimes P_{C(G)}H^d(G).$$
\end{cor}

\subsection{Invariants of restriction to $C(G)$} \label{res subsection}

If $i: C < G$ is a central $p$--elementary abelian of rank $c$, then
\begin{equation*}
H^*(C) \simeq
\begin{cases}
\F_2[x_1,\dots,x_c] & \text{if $p=2$}  \\ \Lambda(x_1,\dots,x_c) \otimes \F_p[y_1,\dots,y_c] & \text{if $p$ is odd},
\end{cases}
\end{equation*}
where $|x_i| = 1$ and $y_i = \beta(x_i)$, and is a Hopf algebra in the usual way.

In \secref{E0 section}, we will see that, after a change of basis for $H^1(C)$, the image of the restriction homomorphism $i^*: H^*(G) \ra H^*(C)$ will be a sub Hopf algebra of $H^*(C)$ of the form
\begin{equation*}
\im(i^*) =
\begin{cases}
\F_2[x_1^{2^{j_1}},\dots,x_c^{2^{j_c}}] & \text{if $p=2$}  \\

\F_p[y_1^{p^{j_1}},\dots,y_b^{p^{j_b}},y_{b+1},\dots,y_c] \otimes \Lambda(x_{b+1},\dots,x_c) & \text{if $p$ is odd},
\end{cases}
\end{equation*}
with the $j_i$ forming a sequence of nonincreasing nonnegative integers\footnote{In the odd prime case, $c-b$ will be the rank of the largest subgroup of $C$ splitting off $G$ as a direct summand.}.

Now suppose that $C = C(G)$. We will say that $G$ has {\em type} $[a_1,\dots,a_c]$ where
\begin{equation*}
(a_1,\dots,a_c) =
\begin{cases}
(2^{j_1}, \dots, 2^{j_c}) & \text{if $p=2$}  \\ (2p^{j_1}, \dots, 2p^{j_b},1, \dots,1) & \text{if $p$ is odd}.
\end{cases}
\end{equation*}

The type of $G$ has the form $[1, \dots, 1]$ if and only if $G = C \times H$, where $Z(H)$ has order prime to $p$.  In all other cases, $a_1 = 2p^k$ for some $k \geq 0$.

Define $e(G)$ and $h(G)$ by
$$ e(G) = \sum_{i=1}^c (a_i-1),$$
and
\begin{equation*}
h(G) \simeq
\begin{cases}
2p^{k-1} & \text{if } a_1 = 2p^k \text{ with } k \geq 1 \\
1 & \text{if } a_1 = 2 \\ 0 & \text{if } a_1=1.
\end{cases}
\end{equation*}

For example, $Q_8 \times \Z/4$ has type $[4,2]$ when $p=2$, so that $e(Q_8 \times \Z/4) = 4$, and $h(Q_8 \times \Z/4)=2$.

\begin{rem} The careful reader will observe that the type of $G$ is just the list of the degrees of the unstable $\A$--algebra generators of $\im(i^*)$, listed in decreasing order, $e(G)$ is the top nonzero degree of the finite dimensional Hopf algebra
$$H^*(C) \otimes_{H^*(G)} \F_p = H^*(C)/(\im(i^{*>0})),$$
and $h(G)$ is the top nonzero degree of the module $\A \cdot H^1(C)$ projected into this Hopf algebra.
\end{rem}

\subsection{$P_{C(G)}H^*(G)$, $d_0(G)$, and $d_1(G)$ when $G$ is $p$--central}

If $G$ is $p$--central with $C=C(G)$, then $P_{C}H^*(G)$ is a finite dimensional unstable algebra.  We identify its top degree and more.

\begin{thm} \label{PC top prim thm} Let $G$ be $p$--central, $C=C(G)$, and $A$ be a Duflot subalgebra of $H^*(G)$.  Then both $P_{C}H^*(G)$ and $Q_AH^*(G)$ are zero in degrees greater than $e(G)$, and one dimensional in degree $e(G)$. Furthermore, $P_{C}H^{e(G)}(G)$ is annihilated by all positive degree elements of the Steenrod algebra, and, if $G$ is a $p$--group, consists of essential\footnote{Recall that $x \in H^*(G)$ is essential if it restricts to zero on all proper subgroups.} cohomology classes.
\end{thm}

The last statement implies the main result of \cite{ak}: $p$--central $p$--groups have nonzero essential cohomology.

This theorem, combined with \corref{Rd cor} and related results, leads to the following calculation.

\begin{thm} \label{PC d0 thm}  Let $G$ be $p$-central.  Then \ $d_0(G) = e(G)$ and $d_1(G) = e(G) + h(G)$.  Furthermore, if $G$ is a finite group with $p$--Sylow subgroup $P$, and $P$ is $p$--central, then $d_0(G) = d_0(P)$ and $d_1(G) = d_1(P)$.
\end{thm}

\begin{cor} \label{PC d0 cor} If $G$ is $p$--central, and $H<G$, then $d_0(H) \leq d_0(G)$ and $d_1(H) \leq d_1(G)$.
\end{cor}

\begin{exs}

\noindent{\bf (a)} As $Q_8$ is $2$--central of type $[4]$, $d_0(Q_8) = 3$ and $d_1(Q_8) = 3+2 = 5$, in agreement with \cite[(II.4.6)]{hls1}.  \\

\noindent{\bf (b)} The hypotheses of $p$--centrality is needed in the last part of the theorem: as observed in \cite[II.4.7]{hls1}, if $G=GL_2(\F_3)$ and $P=SD_{16}$, then $P$ is the $2$--Sylow subgroup of $G$, but $d_0(G)=0 < 2= d_0(P)$, and $d_1(G)=2 < 4 = d_1(P)$.  Similarly, $G$ needs to be $p$--central in the corollary: if $H=Z/4 < D_8=G$, then $d_0(H) = 1 > 0 = d_0(G)$ and $d_1(H) = 2 > 0 = d_1(G)$.  The example $H=\Z/4 < \Z/8 = G$ shows that the inequalities of the corollary can be equalities, even when $H$ is a proper subgroup of a $p$--central $p$--group $G$. \\

\noindent{\bf (c)} The 2-Sylow subgroup $P$ of the simple group $SU(3,4)$ is $2$--central of type $[8,8]$. Thus $d_0(SU(3,4)) = d_0(P) = 14$ and $d_1(SU(3,4)) = d_1(P) = 18$.  Similarly, the 2-Sylow subgroup $Q$ of the simple group $Sz(8)$ is $2$--central of type $[4,4,4]$. Thus $d_0(Sz(8)) = d_0(Q) = 9$ and $d_1(Sz(8)) = d_1(Q) = 11$.  We will see that $P$ and $Q$ have the largest $d_0$ of all $2$--groups of order dividing 64.  For more about the $SU(3,4)$ example, see \secref{ex section}. \\

\noindent{\bf (d)} In \cite{akm}, the authors associate a  2-central Galois group $\mathcal G_{\F}$ to every field $\F$ of characteristic different from $2$ that is not formally real. (They call this the $W$--group of $\F$ because of its connections to the Witt ring $W\F$ \cite{mis}.) From their construction it is easy to deduce that $\mathcal G_{\F}$ has type $[2,\dots,2]$.  Thus $d_0(\mathcal G_{\F}) = r$ and $d_1(\mathcal G_{\F}) = r+1$, where $\mathcal G_{\F}$ has rank $r$.  In particular, the universal group $W$--group $W(n)$ has $d_0(W(n)) = \binom{n+1}{2}$ and $d_1(W(n)) = \binom{n+1}{2}+1$. For more about this example, see \secref{ex section}.
\end{exs}

\subsection{Central essential cohomology}

Our calculation of $d_0(G)$ when $G$ is $p$--central relies on \corref{Rd cor}.  To understand $d_0(G)$ for general $G$, one needs to use the more complicated formula given in \thmref{Rd thm}.  Using some analysis of this already done by us in our companion paper \cite{k4}, we are led to a formula\footnote{Thus far, we have not found an analogous formula for $d_1(G)$.} for $d_0(G)$ that makes use of the following variant of essential cohomology.

We define $Cess^*(G)$, the {\em central essential cohomology} of $G$, to be the kernel of the restriction map
$$ H^*(G) \ra \prod_{C(G)<U} H^*(C_G(U)),$$
where the product is over $p$--elementary abelian subgroups $U$ of $G$ that are strictly bigger than $C(G)$.

This product is over the empty set if $G$ is $p$--central, so the concept is really only interesting when this is not the case.  Furthermore, a theorem of Carlson \cite{carlson depth} implies that $Cess^*(G)$ is nonzero only if the rank of $C(G)$ equals the depth of $H^*(G)$: see \thmref{Cess depth thm} below for the converse.  If this is the case, $Cess^*(G)$ is precisely the {\em depth essential} cohomology of \cite{carlson et al}.

Note that $Cess^*(G)$ has the following structure, compatible in the usual ways: it is an ideal in $H^*(G)$, an unstable module, and an $H^*(C(G))$--comodule.  We have the following general structural results.

\begin{thm} \label{Cess thm}  If $A$ is a Duflot subalgebra of $H^*(G)$, then the following hold. \\

\noindent{\bf (a)}  $Cess^*(G)$ is a finitely generated free $A$--module. \\

\noindent{\bf (b)} The composite $P_C Cess^*(G) \hra Cess^*(G) \era Q_ACess^*(G)$ is monic. \\

\noindent{\bf (c)} The sequence $\displaystyle 0 \ra Q_ACess^*(G) \ra Q_AH^*(G) \ra \prod_{C(G)<U} Q_AH^*(C_G(U))$
is exact.
\end{thm}

Statement (a) implies that $Cess^*(G)$ is a Cohen--MacCauley module, and thus can be viewed as a variant of D.Green's theorem \cite{green3} about the essential cohomology $Ess^*(G)$.  Statement (b) will have application below.  Statement (c) gets us most of the way towards proving the following theorem.

\begin{thm} \label{Cess depth thm} $Cess^*(G) \neq 0$ if and only if the depth of $H^*(G)$ is the rank of $C(G)$.
\end{thm}

The `only if' statement here is just Carlson's theorem.  The `if' statement is a special case of Carlson's Depth Conjecture, and has been previously proved by D.Green under the extra hypothesis that $G$ is a $p$--group \cite{green1}.

\subsection{$d_0(G)$ for general $G$}

Thanks to \thmref{Cess thm}, we can make the following definitions. If $Cess^*(G)$ is nonzero, define $e^{\prime} (G)$ to be the largest $d$ such that $Q_ACess^d(G)$ is nonzero, and $e^{\prime \prime}(G)$ to be the largest $d$ such that $P_{C}Cess^d(G)$ is nonzero.  If $Cess^*(G)=\mathbf 0$, we let $e^{\prime \prime}(G) = e^{\prime}(G)=-1$. Note that \thmref{Cess thm}(b) implies that $e^{\prime \prime}(G) \leq e^{\prime}(G)$.

Using the formula for $\bar R_dH^*(G)$ given in \thmref{Rd thm}, we will prove the following.

\begin{thm} \label{d0 thm}  $ d_0(G) = \max\{ e^{\prime \prime}(C_G(V)) \ | \ V < G \}.$
\end{thm}

\begin{cor} \label{d0 cor}  $ d_0(G) \leq \max\{ e^{\prime}(C_G(V)) \ | \ V < G \}.$
\end{cor}

When computing these maxima, one can restrict to the $p$--elementary abelian groups $V$ which satisfy $V = C(C_G(V))$.\footnote{Since $C_G(V) = C_G(U)$, where $U=C(C_G(V))$.}
The next proposition says that if one has some a priori computation of the depth of $H^*(G)$, one may be able to cut down even more on the $V$'s to be checked.

\begin{prop} \label{depth prop} Assuming that $V = C(C_G(V))$, $Cess^*(C_G(V)) = \mathbf  0$ unless the rank of $V$ is at least equal to the depth of $H^*(G)$.
\end{prop}

This will be proved by combining Carlson's theorem with some $\U$--technology.  This leads to the following generalization of our calculation of $d_0(G)$ for $p$--central $G$.

\begin{cor} \label{cohen--macauley cor}  If $H^*(G)$ is Cohen--Macauley, then
$$d_0(G) = \max\{ e(C_G(V)) \ | \ V < G \text{ is maximal} \}.$$
\end{cor}

This follows from the above, as $C_G(V)$ is $p$--central when $V$ is maximal. Here we have used that, by \thmref{PC top prim thm}, when $G$ is $p$--central,  $e^{\prime \prime}(G) = e^{\prime}(G) = e(G)$.

\begin{conj} \label{e' conjecture} $e^{\prime}(G) < e(G)$ if $G$ is not $p$--central.
\end{conj}

As will be explained in \secref{cess section}, Benson's Strong Regularity Conjecture \cite{b} asserts that certain local cohomology groups $H_{\tilde H^*(G)}^{i,j}(H^*(G))$ vanish.  We connect our conjecture to his.

\begin{prop} \label{conjectures prop} For a fixed finite group $G$, \conjref{e' conjecture} is implied by the Strong Regularity Conjecture.
\end{prop}

Let $G$ have $p$--rank $r$ and $C(G)$ have rank $c$ with $c<r$.  Benson \cite{b} has shown that his conjecture is true if $r-c \leq 2$.  We deduce the next corollary.

\begin{cor} Let $G$ have $p$--rank r, and let $d$ be the depth of $H^*(G)$. If $r-d \leq 2$, then
$$ d_0(G) \leq \max\{ e(C_G(V)) \ | \ V < G \}.$$
\end{cor}

The hypothesis of this corollary applies to all 2--groups of order dividing 64.

\subsection{Calculations when $p=2$}
The Appendix has various tables of values of $d_0(G)$, $d_1(G)$, $e(G)$, $e^{\prime}(G)$, and $e^{\prime \prime}(G)$ for 2--groups of order dividing 64. The tables were compiled by hand using the calculations in \cite{carlson et al}. Their calculations let one immediately determine if $Cess^*(G) \neq 0$, and, when this is the case, one can read off the values of $e(G)$ and $e^{\prime}(G)$, and sometimes $e^{\prime \prime}(G)$.

From our tables, one learns the following about $d_0(G)$ when $p=2$:
\begin{thm} Let $G$ be a finite group with 2--Sylow subgroup $P$ of order dividing 64.  Then $d_0(G) \leq 7$ unless $P$ is isomorphic to either the Sylow subgroup of $SU(3,4)$, in which case $d_0(G) = 14$, or the Sylow subgroup of $Sz(8)$, in which case $d_0(G) = 9$.
\end{thm}

\subsection{Organization of the paper} The rest of the paper is organized as follows.  The nilpotent filtration of $\U$ is reviewed in \secref{U section}, along with basic properties of the functors $\bar R_d$, and the invariants $d_0$ and $d_1$.  Starting from results in \cite{hls1}, in \secref{Rd section} we then deduce the formulae given in \thmref{LF thm} and \thmref{Rd thm}.  In \secref{prim section}, we prove \thmref{prim thm} with a careful analysis of the Lyndon--Serre spectral sequence associated to the group extension $C \ra G \ra G/C$, heavily using that the spectral sequence is a spectral sequence of $H^*(C)$--comodules.  In the $p$--central case, we also input Carlson and Benson's theorem that if $H^*(G)$ is Cohen--Macauley then it is Gorenstein: this leads to proofs of \thmref{PC top prim thm} and \thmref{PC d0 thm} in \secref{pC section}.  Using an analysis of the formula in \thmref{Rd thm} done by us in \cite{k4}, \thmref{d0 thm} is proved in \secref{cess section}, which then continues with our results about $Cess^*(G)$ and the conjectured inequality $e^{\prime}(G) \leq e(G)$.  Though short examples occur throughout, some longer examples that illustrate the general theory make up \secref{ex section}.

\subsection{Acknowledgements}
Early versions of some of these results were presented at a talk in Oberwolfach in September, 2005: see \cite{k3}.  Some of this research was done during a visit to the Cambridge University Pure Mathematics Department in the first half of 2006, and the author thanks DPMMS for its hospitality.  Conversations with Dave Benson about the local cohomology spectral sequence have been enlightening, as have conversations with David Green.

\section{The nilpotent filtration of $\U$} \label{U section}

The nilpotent filtration of $\U$ was introduced in \cite{s1}, and its main properties were developed in \cite{k2,s2,hls1,broto zarati 1}.  Here we collect the results that we need.\footnote{This section necessarily overlaps with the presentation in our recent preprint \cite{k4}.}

\subsection{The definition of $L_d, R_d, \text{ and } \bar R_d$}

For $d\geq 0$, let $\mathcal Nil_d \subset \U$ be the localizing subcategory generated by $d$--fold suspensions of unstable $\A$--modules, i.e. $\mathcal Nil_d$ is the smallest full subcategory containing all $d$--fold suspensions of unstable modules that is closed under extensions and filtered colimits.  Associated to the descending filtration
$$ \dots \subset \mathcal Nil_2 \subset \mathcal Nil_1 \subset \mathcal Nil_0 = \U,$$
there is a natural localization tower for $M \in \U$,
\begin{equation*}
\xymatrix{
&&& \vdots \ar[d] \\
&&& L_2M \ar[d] \\
&&& L_1M \ar[d] \\
M \ar[rrr]^{\lambda_0}  \ar[urrr]^{\lambda_1} \ar[uurrr]^{\lambda_2} &&& L_0M,
}
\end{equation*}
where $L_d: \U \ra \U$ is localization away from $\mathcal Nil_{d+1}$.\footnote{What we are calling $L_d$ here was called $L_{d+1}$ in \cite{hls1}.}  The natural transformation $\lambda_d:M \ra L_dM$ is characterized by the following properties: \\

\noindent{\bf (a)} \  $L_dM$ is $\mathcal Nil_{d+1}$--closed, i.e. $\Ext^s_{\U}(N,L_dM) = 0$ for $s=0,1$ and $N \in \mathcal Nil_{d+1}$, \\

\noindent{\bf (b)} \ $\lambda_d$ is a $\mathcal Nil_{d+1}$--isomorphism, i.e. $\ker \ \lambda_d$ and $\coker \ \lambda_d$ are both in $\mathcal Nil_{d+1}$. \\

A module $M \in \U$ admits a natural filtration
$$ \dots \subseteq nil_2M \subseteq nil_1M \subseteq nil_0 M= M,$$
where $nil_dM$ is the largest submodule in $\mathcal Nil_d$.  For $d>0$, $nil_dM = \ker \ \lambda_{d-1}$.

An unstable module $M$ is called {\em reduced} if $nil_1M=0$.  As observed in \cite[Prop.2.2]{k2},
$ nil_dM/nil_{d+1}M = \Sigma^d R_dM$,
where $R_dM$ is a reduced unstable module. (See also \cite[Lemma 6.1.4]{s2}.) Then $\bar R_dM$ is defined to be the $\mathcal Nil_1$--closure of $R_dM$.  Thus $R_dM \subseteq L_0R_dM = \bar R_dM$.

We have the following useful alternative definition of $\bar R_dM$.  (Compare with \cite[I(3.8.1)]{hls1}.)

\begin{prop} \label{loc prop}   There is a natural isomorphism
$$\Sigma^d \bar R_dM \simeq \ker \ \{ L_dM \ra L_{d-1}M \}.$$
\end{prop}

The functors $L_d$ and $L_{d-1}$ are left exact, as they are localizations, and thus we conclude

\begin{cor}  $\bar R_d: \U \ra \U$ is left exact.
\end{cor}

\begin{proof}[Proof of \propref{loc prop}]  Let $c_{d+1}M = \coker \ \{ \lambda_d: M \ra L_dM \}$. Then $c_{d+1}M \in \mathcal Nil_{d+1}$, and there is an exact sequence
$$ 0 \ra nil_{d+1}M \ra M \ra L_dM \ra c_{d+1}M \ra 0.$$
Diagram chasing then shows that there is a natural short exact sequence
$$ 0 \ra nil_dM/nil_{d+1}M \ra \ker \ \{ L_dM \ra L_{d-1}M \} \ra \ker \ \{ c_{d+1}M \ra c_dM \} \ra 0.$$
As the middle module here is $\mathcal Nil_{d+1}$--closed, and the right module is in $\mathcal Nil_{d+1}$, we see that the left map identifies with $\lambda_d$.  Recalling that $ nil_dM/nil_{d+1}M = \Sigma^d R_dM$, this says that there is a natural isomorphism
$$ L_d(\Sigma^d R_dM) \simeq \ker \ \{ L_dM \ra L_{d-1}M \}.$$
The proof of the proposition is then completed by observing that $L_d(\Sigma^d R_dM) \simeq \Sigma^d \bar R_dM$, a consequence of the next proposition.
\end{proof}

\begin{prop} \label{susp loc prop} There is a natural isomorphism $L_{c+d}(\Sigma^dM) \simeq \Sigma^d L_cM$, for all $M \in \U$.
\end{prop}

\begin{proof}  We need to check that the map $\Sigma^d \lambda_c: \Sigma^dM \ra \Sigma^d L_cM$ satisfies the two properties characterizing localization away from $\mathcal Nil_{c+d+1}$.

That $ker(\Sigma^d \lambda_c)$ and $\coker(\Sigma^d \lambda_c)$ are both in $\mathcal Nil_{c+d+1}$ is clear, as $\ker(\lambda_c)$ and $\coker(\lambda_c)$ are both $\mathcal Nil_{c+1}$, and the $d$--fold suspension of a module in $\mathcal Nil_{c+1}$ will be in $\mathcal Nil_{c+d+1}$.

To see that the range of $\Sigma^d \lambda_c$  is $\mathcal Nil_{c+d+1}$--closed, we check that if $M \in \U$ is $\mathcal Nil_{c+1}$--closed then $\Sigma^d M$ is $\mathcal Nil_{c+d+1}$--closed.  This follows from the following characterization of $\mathcal Nil_{c+1}$--closed modules: $M \in \U$ is $\mathcal Nil_{c+1}$--closed if and only if it fits into an exact sequence of the form
$$ 0 \ra M \ra \prod_{\alpha} H^*(V_{\alpha}) \otimes M_{\alpha}  \ra \prod_{\beta} H^*(W_{\beta}) \otimes N_{\beta},$$
with all the modules $M_{\alpha}$ and $N_{\beta}$ concentrated in degrees between $0$ and $c$.  See \cite[Prop.1.15]{broto zarati 2}.
\end{proof}

\subsection{Further properties of $L_d, R_d, \text{ and } \bar R_d$}

We need to recall some notation and terminology. If $V$ is an elementary $p$--group, $T_V: \U \ra \U$ is defined to be the left adjoint to $H^*(V) \otimes \text{\underline{\hspace{.1in} }}$, as famously studied by Lannes \cite{L1,L3}.  Given a Noetherian unstable algebra $K \in \K$, $K_{f.g.}-\U$ is defined to be the category studied in \cite[I.4]{hls1} whose objects are finitely generated $K$--modules $M$ whose $K$--module structure map $K \otimes M \ra M$ is in $\U$, and morphisms are $K$--module maps in $\U$.

\begin{prop} \label{Ld properties} The functor $L_d: \U \ra \U$ satisfies the following properties. \\

\noindent{\bf (a)} There are natural isomorphisms $ L_0(M \otimes N) \simeq L_0M \otimes L_0N$. \\

\noindent{\bf (b)} There are natural isomorphisms $T_V L_dM \simeq L_dT_V M$. \\

\noindent{\bf (c)} If $K \in \K$, then $L_dK \in \K$, and $K \ra L_dK$ is a map of unstable algebras.  If $K$ is also Noetherian, and $M \in K_{f.g.}-\U$, then $L_dK \in K_{f.g.}-\U$, and thus is Noetherian, and $L_dM \in L_dK_{f.g.}-\U$.
\end{prop}

Property (b) can be deduced from properties of $T_V$ as follows.  First, to see that $T_VL_dM$ is $\mathcal Nil_{d+1}$--closed, we compute, for $s=0,1$ and $N \in \mathcal Nil_{d+1}$:
$$ \Ext^s_{\U}(N,T_VL_dM) = \Ext^s_{\U}(H^*(V) \otimes N,L_dM) = 0, $$
since $H^*(V) \otimes N$ will be in  $\mathcal Nil_{d+1}$ if $N$ is.  Second,
$ T_V \lambda_d: T_VM \ra T_VL_dM$
is a $\mathcal Nil_{d+1}$--isomorphism, as the kernel and cokernel are in $\mathcal Nil_{d+1}$, since $T_V$ is exact and sends $\mathcal Nil_{d+1}$ to itself.

See \cite[I.4]{hls1} and \cite{broto zarati 1} for more detail about properties (a) and (c).

\begin{prop} \label{Rd properties} The functors $R_d: \U \ra \U$ satisfy the following properties. \\

\noindent{\bf (a)} There are a natural isomorphisms $ R_*(M \otimes N) \simeq R_*M \otimes R_*N$ of graded objects in $\U$.  \\

\noindent{\bf (b)} There are natural isomorphism $T_V R_dM \simeq R_dT_V M$.  \\

\noindent{\bf (c)} If $K \in \K$, then $R_0K \in \K$, and $K \ra R_0K$ is a map of unstable algebras.  If $K $ is also Noetherian, and $M \in K_{f.g.}-\U$, then $R_0K$ is also a Noetherian unstable algebra, and $R_dM \in R_0K_{f.g.}-\U$, for all $d$.
\end{prop}

For the first two properties, see \cite[\S 3]{k2}, and the last follows easily from the first.

\begin{prop} \label{bar Rd properties} The functors $\bar R_d: \U \ra \U$ satisfy the following properties. \\

\noindent{\bf (a)} There are natural isomorphisms $\bar R_*(M \otimes N) \simeq \bar R_*M \otimes \bar R_*N$ of graded objects in $\U$. \\

\noindent{\bf (b)} There are natural isomorphisms $T_V \bar R_dM \simeq \bar R_dT_V M$.  \\

\noindent{\bf (c)} If $K \in \K$, then $\bar R_0K \in \K$, and $K \ra \bar R_0K$ is a map of unstable algebras.  If $K$ is also Noetherian, and $M \in K_{f.g.}-\U$, then $\bar R_0K$ is also a Noetherian unstable algebra, and $\bar R_dM \in \bar R_0K_{f.g.}-\U$, for all $d$.
\end{prop}

This, of course, follows from the previous two propositions.

A Noetherian unstable algebra $K$ has a finite Krull dimension $\dim K$.  We have an addendum to \propref{Ld properties}.

\begin{prop}[{\cite[Prop.4.10]{k4}}] \label{Krull dim prop} If an unstable algebra $K$ is Noetherian, then $\dim K = \dim L_0K$.
\end{prop}

Another special property of $L_0$ that we will need goes as follows.

\begin{prop}[{\cite[Lem.4.3.3]{L1}}] \label{L0 prop}  Let $f:M \ra N$ be a map in $\K$.  Then $$L_0f:L_0M \ra L_0N$$ is an isomorphism if and only if, for all $p$--elementary abelian groups $V$, the induced map
$$ f^*: \Hom_{\K}(N,H^*(V)) \ra \Hom_{\K}(M,H^*(V))$$
is a bijection.
\end{prop}

As in the introduction,  given $M \in \U$, $M_{LF}$ denotes the submodule of locally finite elements: $x \in M$ such that $\A x \subseteq M$ is finite.

\begin{prop} \label{LF and Rd prop} There is a natural isomorphism $(R_d M)^0 = (\bar R_dM)^0 \simeq (M_{LF})^d$.
\end{prop}

See \cite[\S 3]{k2} for a proof.

Finally,  Henn \cite{henn} proved the following important finiteness result.

\begin{prop} \label{finite prop 2} Let $K \in \K$ be Noetherian, and $M \in K_{f.g.}-\U$. Then the $M$ is $\mathcal Nil_{d}$--local for $d >> 0$.  In particular, the nilpotent filtration of $M$ has finite length.
\end{prop}

\subsection{Properties of $d_0M$ and $d_1M$}

The authors of \cite{hls1} define $d_0M$ and $d_1M$ as follows.

\begin{defn}  Let $M$ be an unstable module. \\

\noindent{\bf (a)} \ Let $d_0M$ be the smallest $d$ such that $\lambda_d$ is monic, or $\infty$ if no such $d$ exists.  Equivalently, $d_0M$ is the smallest $d$ such that $\Hom_{\U}(N,M) = 0$ for all $N \in \mathcal Nil_{d+1}$, or the smallest $d$ such that $nil_{d+1}M = 0$. If $M$ is nonzero, $d_0M$ is also the largest $d$ such that $R_dM$ is nonzero, or the largest $d$ such that $\bar R_dM$ is nonzero. \\

\noindent{\bf (b)} \ Let $d_1M$ be the smallest $d$ such that $\lambda_d$ is an isomorphism, or $\infty$ if no such $d$ exists.  Equivalently, $d_1M$ is the smallest $d$ such that $\Ext^s_{\U}(N,M) = 0$ for $s=0,1$ and all $N \in \mathcal Nil_{d+1}$. \\
\end{defn}

As fundamental examples, we have that $d_0H^*(V) = d_1H^*(V) = 0$ for all elementary abelian $p$--groups $V$.

\begin{prop} \label{d1 prop 1} Let $M$ and $N$ be unstable modules. \\

\noindent{\bf (a)} \ For $s=0,1$, $d_s(M \oplus N) = \max \{ d_sM,d_sN \}$. \\

\noindent{\bf (b)} \ If $M$ and $N$ are nonzero, $d_0(M \otimes N) = d_0M + d_0N$ and $d_1(M \otimes N) = \max \{ d_1M + d_0N,d_0M +d_1N \}$. \\

\noindent{\bf (c)} \ For $s=0,1$, $d_sT_VM = d_sM$. \\

\noindent{\bf (d)} \ If $M$ is nonzero, for $s=0,1$, $d_s(\Sigma^n M) = d_sM + n$. \\
\end{prop}

For properties (a) and (b) see \cite[Prop.I.3.6]{hls1}.  Using the exactness of $T_V$, property (c) follows from \propref{Ld properties}(b).  Property (d) follows from \propref{susp loc prop}.

\begin{prop}  \label{d1 prop 2} Let $0 \ra M_1 \ra M_2 \ra M_3 \ra 0$ be a short exact sequence in $\U$. \\

\noindent{\bf (a)} For $s=0,1$, $d_sM_2 \leq \max \{d_sM_1,d_sM_3 \}$.  Furthermore, if $d_sM_3 < d_sM_1$, then $d_sM_2 = d_sM_1$. \\

\noindent{\bf (b)} $d_0M_1 \leq d_0M_2$ and $d_1M_1 \leq \max \{d_1M_2, d_0M_3 \}$.  Furthermore, if $d_1M_2 < d_0M_3$, then $d_1M_1=d_0M_3$.  \\
\end{prop}

This is proved with straightforward use of the long exact $\Ext^*$ sequence associated to a short exact sequence.  Compare with \cite[Prop.I.3.6]{hls1}.  \\

\begin{cor} \label{d1 cor} If $M \in \U$ is reduced, then $d_1M = d_0(L_0M/M)$. \\
\end{cor}

This follows by applying \propref{d1 prop 2}(b) to $0 \ra M \ra L_0M \ra L_0M/M \ra 0$.  \\

\subsection{Basic properties of $d_0(G)$ and $d_1(G)$}

 By abuse of notation, if $G$ is a finite group, for $s=0,1$, we write $d_s(G)$ for $d_sH^*(G)$. For example, $d_0(V)=d_1(V)=0$ for all elementary abelian $p$--groups $V$.

 The properties of $d_0M$ and $d_1M$ presented above have the following immediate consequences for $d_0(G)$ and $d_1(G)$

 \begin{prop} Let $G$ and $H$ be finite groups. \\

 \noindent{\bf (a)} $d_0(G \times H) = d_0(G) + d_0(H)$. \\

 \noindent{\bf (b)} $d_1(G \times H) = \max \{ d_1(G)+ d_0(H),d_0(G) +d_1(H) \}$. \\

 \noindent{\bf (c)} If $P$ is a $p$--Sylow subgroup of $G$, then $d_s(G) \leq d_s(P)$ for $s=0,1$. \\

 \noindent{\bf (d)} If $V$ is a $p$--elementary abelian subgroup of $G$, then  $d_s(C_G(V)) \leq d_s(G)$ for $s=0,1$.
 \end{prop}

 Properties (a) and (b) follow from \propref{d1 prop 1} (b).  As the unstable module $H^*(G)$ is a direct summand of $H^*(P)$ if $P$ is a $p$--Sylow subgroup, property (c) follows from \propref{d1 prop 1} (a).  Similarly property (d) follows from \propref{d1 prop 1} (c), as $H^*(C_G(V))$ is a direct summand of $T_VH^*(G)$ \cite{lannes 2}\footnote{It is unfortunate that this much referenced elegant 1986 preprint has never been published.}.

\section{Formulae for $H^*(G)_{LF}$ and $\bar R_d(H^*(G))$} \label{Rd section} In this section we prove the formulae for $H^*(G)_{LF}$ and $\bar R_dH^*(G)$ given in \secref{main results section}.

\subsection{A formula for $L_dH^*(G)$}

The starting point for all of these are the following constructions.  Given a morphism $\alpha: V_1 \ra V_2$ in $\A(G)$, there are maps
$$ \alpha^*: H^*(V_2) \ra H^*(V_1),$$
$$ \alpha_*: H^*(C_G(V_1)) \ra H^*(C_G(V_2)), \text{ and } $$
$$ m_{\alpha}^*: H^*(C_G(V_2)) \ra H^*(V_1)\otimes H^*(C_G(V_2)).$$
Here $\alpha_*$ is induced by conjugation by $g^{-1}$ where $g \in G$ is any element\footnote{This is well defined as any two choices will differ by an element of $C_G(V_1)$, and so will agree on cohomology.} chosen so that conjugation by $g$ induces $\alpha$, and
$$ m_{\alpha}: V_1 \times C_G(V_2) \ra C_G(V_2)$$
is the homomorphism sending $(x,y)$ to $\alpha(x)y$.
We also let $m_V: V \times V \ra V$ denote multiplication in an elementary abelian group $V$.

To state one of the formulae from \cite{hls1}, we recall two other bits of notation from \secref{main results section}.  Given an unstable module $M$, we let $M^{\leq d}$ denote $M$ modulo degrees greater than $d$. Given a category $\mathcal C$, we let $\mathcal C^{\#}$ denote the associated twisted arrow category:
the objects of $\mathcal C^{\#}$ are the morphisms of $\mathcal C$, and a morphism $\alpha \rightsquigarrow \beta$ from $\alpha: A_1 \ra A_2$ to $\beta: B_1 \ra B_2$ is a commutative diagram in $\mathcal C$
\begin{equation*}
\xymatrix{
A_1  \ar[d] \ar[r]^{\alpha} & A_2   \\
B_1 \ar[r]^{\beta} & B_2. \ar[u] }
\end{equation*}

\cite[Formula I(5.5.1)]{hls1} now reads

\begin{thm} \label{Ld formula} The homomorphisms $V_1 \times C_G(V_2) \xra{m_{\alpha}} C_G(V_2) \subset G$ induce an isomorphism of unstable algebras from $L_d H^*(G)$ to
$$ \lim_{V_1 \xra{\alpha} V_2} \text{Eq } \{\xymatrix{H^*(V_1)\otimes H^{\leq d}(C_G(V_2)) \ar@<.5ex>[r]^-{\mu(\alpha)} \ar@<-.5ex>[r]_-{\nu(\alpha)} & H^*(V_1) \otimes (H^*(V_1) \otimes H^*(C_G(V_2)))^{\leq d}}\},$$
where $\mu(\alpha)$ is induced by $1 \otimes m_{\alpha}^*$, $\nu(\alpha)$ is induced by $m_{V_1}^*\otimes 1$, and the limit is over $\A(G)^{\#}$.
\end{thm}

\subsection{A formula for $\bar R_dH^*(G)$}

Recall our notation from \secref{main results section}: if $W$ is a central elementary abelian $p$--subgroup of $Q$, then $P_WH^*(Q)$ denotes the algebra of primitives in the $H^*(W)$--comodule $H^*(Q)$.

\begin{prop} \label{Rd prop} As unstable modules, $\bar R_d H^*(G)$ is naturally isomorphic to
$$ \lim_{V_1 \xra{\alpha} V_2} H^*(V_1) \otimes P_{\alpha(V_1)}H^d(C_G(V_2)),$$
where the limit is over $\A(G)^{\#}$.
\end{prop}
\begin{proof}  Recall that $\Sigma^d \bar R_d M$ is the kernel of $L_dM \ra L_{d-1}M$. As kernels commute with limits and equalizers, it follows from the previous theorem that $\bar R_d H^*(G)$ is naturally isomorphic to
\begin{equation*}  \lim_{V_1 \xra{\alpha} V_2} \text{Eq } \{ \xymatrix{H^*(V_1)\otimes H^d(C_G(V_2)) \ar@<.5ex>[r]^-{\mu(\alpha)} \ar@<-.5ex>[r]_-{\nu(\alpha)} & H^*(V_1) \otimes (H^*(V_1) \otimes H^*(C_G(V_2)))^{d}}\},
\end{equation*}
where $\mu(\alpha)$ is induced by $1 \otimes m_{\alpha}^*$ and $\nu(\alpha)$ is induced by $m_{V_1}^*\otimes 1$.  But now we observe that the equalizer in this formula is precisely $H^*(V_1)\otimes P_{\alpha(V_1)}H^d(C_G(V_2))$.  For $\nu(\alpha)$ is the composite
\begin{equation*}
\begin{split}
H^*(V_1)\otimes H^d(C_G(V_2)) &
\xra{m_{V_1}^*\otimes 1} H^*(V_1) \otimes H^*(V_1) \otimes H^d(C_G(V_2)) \\
  & \xra{\text{truncate}}  H^*(V_1) \otimes H^0(V_1) \otimes H^d(C_G(V_2)),
\end{split}
\end{equation*}
and this identifies with
$$ H^*(V_1)\otimes H^d(C_G(V_2)) \xra{1 \otimes \pi^*}   H^*(V_1) \otimes (H^*(V_1) \otimes H^*(C_G(V_2)))^d,$$
where $\pi: V_1 \times C_G(V_2) \ra C_G(V_2)$ is the projection.
\end{proof}

\subsection{A formula for $H^*(G)_{LF}$}

\begin{prop} \label{LF prop} As unstable algebras, $H^*(G)_{LF}$ is naturally isomorphic to
$$ \lim_{V_1 \xra{\alpha} V_2} P_{\alpha(V_1)}H^*(C_G(V_2)),$$
where the limit is over $\A(G)^{\#}$.
\end{prop}
\begin{proof}  As there are no nonzero locally finite elements in $\tilde H^*(V_1) \otimes H^*(C_G(V_2))$, the composite $H^*(G)_{LF} \subset H^*(G) \ra H^*(C_G(V_2))$ has image in $P_{\alpha(V_1)}H^*(C_G(V_2))$ for any $\alpha: V_1 \ra V_2$ in $\A(G)$.  Thus one gets a natural map of unstable algebras
$$ H^*(G)_{LF} \ra \lim_{V_1 \xra{\alpha} V_2} P_{\alpha(V_1)}H^*(C_G(V_2)).$$
That this is an isomorphism follows from \propref{Rd prop}, recalling that \propref{LF and Rd prop} said that there is a natural isomorphism $(\bar R_d M)^0 \simeq (M_{LF})^d$.
\end{proof}

\subsection{Replacing $\A(G)$ with $\A_C(G)$}

Recall that $C(G)$ denotes the maximal central $p$--elementary abelian subgroup of $G$, and $\A_C(G)$ denotes the full subcategory of $\A(G)$ with objects $C(G) \leq V < G$.

\begin{thm} \label{AC theorem} One can take the limit over $\A_C(G)^{\#}$, rather than $\A(G)^{\#}$ in \thmref{Ld formula}, \propref{Rd prop}, and \propref{LF prop}.
\end{thm}

This will follow quite formally from the following simple observations.  Let $C=C(G)$. Given $V < G$, let $CV < G$ be the subgroup generated by $C$ and $V$.  This induces an evident functor $C: \A(P) \ra \A_C(G)$. Furthermore, the natural inclusion $V \ra CV$ induces an identification $C_G(CV) = C_G(V)$.

Given $\alpha: V_1 \ra V_2$, let $\alpha_C: V_1 \ra CV_2$ be the evident map, and then let
$$\xymatrix{\alpha & \alpha_C \ar[l]_-{f_{\alpha}} \ar[r]^-{g_{\alpha}} & C\alpha, }$$
morphisms in $\A(G)^{\#}$, correspond to the diagram in $\A(G)$
\begin{equation*}
\xymatrix{
V_1 \ar[d]^{\alpha} \ar@{=}[r] & V_1 \ar[d]^{\alpha_C} \ar[r] & CV_1 \ar[d]^{C\alpha} \\
V_2 \ar[r] & CV_2 \ar@{=}[r] & CV_2.}
\end{equation*}

\begin{lem}  \label{AC lemma} Let $F: \A(G)^{\#} \ra \F_p\text{-vector spaces}$ be a contravariant functor such that for all $\alpha:V_1 \ra V_2$, $F(f_{\alpha}): F(\alpha) \ra F(\alpha_C)$ is an isomorphism.  Then the natural map
$$ \Psi: \lim_{\alpha \in \A(G)^{\#}}F(\alpha) \ra \lim_{\alpha \in \A_C(G)^{\#}}F(\alpha)$$
is an isomorphism.
\end{lem}

Note that both $F(V_1 \xra{\alpha} V_2) = H^*(V_1)$ and $F(V_1 \xra{\alpha} V_2) = H^*(C_G(V_2))$ satisfy the hypothesis of the lemma.  \thmref{AC theorem} then follows from the lemma, as the relevant $F$'s are built from these two examples by constructions that preserve isomorphisms.

\begin{proof}[Proof of \lemref{AC lemma}]  We define $\displaystyle \Phi: \lim_{\alpha \in \A^C(G)^{\#}}F(\alpha) \ra \lim_{\alpha \in \A(P)^{\#}}F(\alpha)$, an inverse to $\Psi$, as follows.  Given $\displaystyle x = (x_{\beta}) \in \lim_{\beta \in \A^C(G)^{\#}}F(\beta)$, let $ \Phi(x) = (\Phi(x)_{\alpha}) \in \prod_{\alpha \in \A(G)^{\#}}F(\alpha)$, where $\Phi(x)_{\alpha} = F(f_{\alpha})^{-1}F(g_{\alpha})(x_{C\alpha})$.  One then checks that $\displaystyle \Phi(x) \in \lim_{\alpha \in \A(G)^{\#}}F(\alpha)$, $\Psi \circ \Phi = 1$, and $\Phi \circ \Psi = 1$.
\end{proof}

\subsection{Rewriting the formulae}

If $\mathcal C$ is a small category, and
$$F: {\mathcal C}^{\#} \ra \F_p\text{-vector spaces}$$ is a contravariant functor, there is a canonical isomorphism
$$ \lim_{{\mathcal C}^{\#}} F = \text{Eq} \left\{ \prod_{C \in ob \mathcal C} F(1_C)
\begin{array}{c} \mu \\[-.08in] \longrightarrow \\[-.1in] \longrightarrow \\[-.1in] \nu
\end{array}
\prod_{\alpha \in mor \mathcal C} F(\alpha) \right\},
$$
where, given $\alpha: C_1 \ra C_2$, the $\alpha$-component of $\mu$ and $\nu$ are induced by applying $F$ to the canonical morphisms in ${\mathcal C}^{\#}$ from $\alpha$ to $1_{C_1}$ and $1_{C_2}$ respectively.

Thus, for example, $\bar R_dH^*(G)$ will be naturally isomorphic to
$$ \text{Eq} \left\{ \prod_{V} H^*(V) \otimes P_VH^d(C_G(V))
\begin{array}{c} \mu \\[-.08in] \longrightarrow \\[-.1in] \longrightarrow \\[-.1in] \nu
\end{array}
\prod_{\alpha:V_1 \ra V_2} H^*(V_1) \otimes P_{\alpha(V_1)}H^d(C_G(V_2)) \right\},
$$where $\mu$ and $\nu$ are induced by
$$ 1 \otimes \alpha_*: H^*(V_1) \otimes P_{V_1}H^d(C_G(V_1)) \ra H^*(V_1) \otimes P_{\alpha(V_1)}H^d(C_G(V_2))$$
and
$$ \alpha^* \otimes i: H^*(V_2) \otimes P_{V_2}H^d(C_G(V_2)) \ra H^*(V_1) \otimes P_{\alpha(V_1)}H^d(C_G(V_2))$$
for each $\alpha: V_1 \ra V_2$. ($i$ is the evident inclusion.)

Morphisms in $\A(G)$ factor as inclusions followed by isomorphisms induced by the inner automorphism group $Inn(G)$, so this last formula rewrites as follows.

\begin{prop} \label{good Rd formula prop} $\bar R_dH^*(G)$ is naturally isomorphic to
$$ \text{Eq} \left\{ \left[\prod_{V} H^*(V) \otimes P_VH^d(C_G(V))\right]^{Inn(G)}
\begin{array}{c} \mu \\[-.08in] \rightarrow \\[-.1in] \rightarrow \\[-.1in] \nu
\end{array}
\prod_{V_1 < V_2} H^*(V_1) \otimes P_{V_1}H^d(C_G(V_2)) \right\},
$$
where $\mu$ and $\nu$ are induced by
$$ 1 \otimes \eta_*: H^*(V_1) \otimes P_{V_1}H^d(C_G(V_1)) \ra H^*(V_1) \otimes P_{V_1}H^d(C_G(V_2))$$
and
$$ \eta^* \otimes i: H^*(V_2) \otimes P_{V_2}H^d(C_G(V_2)) \ra H^*(V_1) \otimes P_{V_1}H^d(C_G(V_2))$$
for each inclusion $\eta: V_1 < V_2$ in $\A_C(G)$.
\end{prop}

Otherwise said, $$ x = (x_V) \in \left[\prod_V H^*(V) \otimes P_VH^d(C_G(V))\right]^{Inn(G)}$$ is in $\bar R_dH^*(G)$ exactly when the components are related by $$(1 \otimes \eta_*)(x_{V_1}) = (\eta^* \otimes i)(x_{V_2}),$$ for each inclusion $\eta: V_1 < V_2$ in $\A_C(G)$.

Similarly, we have

\begin{prop} \label{good LF formula prop} $H^*(G)_{LF}$ is naturally isomorphic to
$$ \text{Eq} \left\{ \left[\prod_{V} P_VH^*(C_G(V))\right]^{Inn(G)}
\begin{array}{c} \mu \\[-.08in] \longrightarrow \\[-.1in] \longrightarrow \\[-.1in] \nu
\end{array}
\prod_{V_1 < V_2} P_{V_1}H^*(C_G(V_2)) \right\},
$$
where $\mu$ and $\nu$ are induced by
$$ \eta_*: P_{V_1}H^*(C_G(V_1)) \ra P_{V_1}H^*(C_G(V_2))$$
and
$$ i: P_{V_2}H^*(C_G(V_2)) \subseteq P_{V_1}H^*(C_G(V_2))$$
for each inclusion $\eta: V_1 < V_2$ in $\A_C(G)$.
\end{prop}

\section{The cohomology of central extensions} \label{prim section}

Let $C$ be a central $p$--elementary abelian subgroup of a finite group $G$, and let $Q = G/C$.  This is the first of two sections in which we study the rich structure of the Lyndon--Hochschild--Serre spectral sequence $\{E_r^{*,*}(G,C) \}$ associated to the central extension:
$$ C \xra{i} G \xra{q} Q.$$
Enroute, we will prove \thmref{prim thm}.

To begin, we recall the following standard facts. The spectral sequence is a spectral sequence of differential graded algebras, converging to $H^*(G)$, and with $E_2^{*,*} = H^*(Q) \otimes H^*(C)$.  Furthermore, $E_r^{*,*} = E_{\infty}^{*,*}$ for $r >>0$ \cite{evens}.

Recall that the extension corresponds to an element $\tau \in H^2(Q;C)$, or equivalently  a homomorphism $\tau:H_2(Q) \ra C$. Under the identification $C^{\#} = H^1(C)$, its dual $\tau^{\#}: C^{\#} \ra H^2(Q)$ corresponds to $ d_2: E_2^{0,1} \ra E_2^{2,0}$.

\subsection{$H^*(C)$--comodule structure of the spectral sequence}

As $C$ is central, multiplication $m: C \times G \ra G$ is a group homomorphism. The induced  algebra map
$$ m^*: H^*(G) \ra H^*(C) \otimes H^*(G), $$
makes $H^*(G)$ into a $H^*(C)$--comodule.  The restriction $i^*: H^*(G) \ra H^*(C)$ is both an algebra and comodule map, and it follows that  $E_{\infty}^{0,*} = \im(i^*)$ is a subHopf algebra of $E_{2}^{0,*} = H^*(C)$.

One can strengthen these last observations to statements about the whole spectral sequence.  A good functorial model for $BG$, say the reduced bar construction, shows that $BC$ is an abelian topological group, $BG$ is a $BC$--space equipped with proper free action via $Bm: BC \times BG \ra BG$, and $BG \ra BQ$ is the associated principal $BC$--bundle.  The Serre spectral sequence arises from the pullback to $BG$ of the skeletal filtration of $BQ$.  This will be a filtration of $BG$ by $BC$--subspaces, and we conclude the following.

\begin{lem}  For all $k$ and $r$, $E_r^{k,*}$ is an $H^*(C)$--comodule, such that the maps
$$ d_r: E_r^{k,*} \ra E_r^{k+r,*}$$
and
$$  E_r^{i,*} \otimes E_{r}^{j,*} \ra E_r^{i+j,*}$$
are maps of $H^*(C)$--comodules.   In particular, $E_{r}^{0,*}$ is a subHopf algebra of $E_{2}^{0,*} = H^*(C)$.
\end{lem}

\subsection{A handy Hopf algebra lemma}

We now digress to state and prove a handy statement about (connected graded) Hopf algebras that we can apply to the situation of the previous subsection.

We need some notation. Let $H$ be a graded connected Hopf algebra over a field $\F$. There is a canonical splitting of vector spaces $H = \F \oplus I(H)$, where $I(H)$ is the augmentation ideal.   If $M$ is a right $H$--module, let the module of indecomposables be defined by $Q_HM = M \otimes_H \F = M/MI(H)$.  Dually, if $M$ is a right $H$--comodule, let the module of primitives be defined by
\begin{equation*}
\begin{split}
P_HM & = \text{Eq } \{ M
\begin{array}{c} \Delta \\[-.08in] \longrightarrow \\[-.1in] \longrightarrow \\[-.1in] i
\end{array}
M \otimes H \} \\
& = \ker \{ \bar{\Delta}: M \ra M \otimes I(H) \},
\end{split}
\end{equation*}
where $\Delta: M \ra M \otimes H$ is the comodule structure, $i$ is the inclusion induced by the unit $\F \ra H$, and $\bar{\Delta}$ is the composite $M \xra{\Delta} M \otimes H \ra M \otimes I(H)$.

\begin{lem} \label{handy lemma}  Let $K$ be a subHopf algebra of a Hopf algebra H.  Suppose $M$ is simultaneously an $H$--comodule and $K$--module such that the $K$--module structure map $ M \otimes K \ra M$ is a map of $H$--comodules.  Then \\

\noindent{\bf (a)} $M$ is a free $K$--module, and \\

\noindent{\bf (b)} the composite $P_HM \hra M \era Q_KM$ is monic.
\end{lem}

\begin{rem} To put this in perspective, the lemma has long been known if $K=H$, and, in this case, $P_HM \simeq Q_HM$ \cite[Thm.4.1.1]{sweedler}.  Our proof is very similar to the proofs of Proposition 1.7 and Theorem 4.4 of Milnor and Moore's classic paper \cite{milnor moore}. Compare also to Green's lemma \cite[Lem2.1]{green3}.
\end{rem}

Before proving the lemma, we note the following consequence.  Given $H$ and $K$ as in the lemma, let $K-H-\mathcal Mod$ be the category of $M$ as in the lemma: an object is a vector space $M$ that is simultaneously an $H$--comodule and $K$--module such that the $K$--module structure map $ M \otimes K \ra M$ is a map of $H$--comodules.  Morphisms are linear maps that are both $K$--module and $H$--comodule maps. $K-H-\mathcal Mod$ is an abelian category in the obvious way.

\begin{cor} \label{K-H-cor}  {\bf (a)} Every short exact sequence $0 \ra M_1 \ra M \ra M_2 \ra 0$ in $K-H-\mathcal Mod$ is split as a sequence of $K$--modules. \\

\noindent{\bf (b)} The functor sending $M$ to $Q_K(M)$ is exact on $K-H-\mathcal Mod$.
\end{cor}

\begin{proof}[Proof of \lemref{handy lemma}]  Choose a section $s: Q_KM \ra M$ of the quotient $\pi: M \ra Q_KM$, and let
$ m_s: Q_KM \otimes I(K) \ra MI(K)$
be the epimorphism given by $m_s(x,k) = s(x)k$.  Statement (a) is asserting that $m_s$ is an isomorphism.

Let $\Delta_K: MI(K) \ra M \otimes I(H)$ be the composite
$$ MI(K) \subset M \xra{\bar{\Delta}}  M \otimes I(H),$$
Statement (b) asserts that $P_HM \cap MI(K) = \{0\}$, i.e. that $\Delta_K$ is monic.

Thus both statements will follow from the following claim:
$$ \Delta_K \circ m_s: Q_KM \otimes I(K) \ra M \otimes I(H)$$
is monic.

To prove this claim, let $F_nM$ be the $K$--submodule of $M$ generated by elements of degree up to $n$.  Given $x \in (Q_KM)^n$, and $k \in I(K)$, let $\Delta(s(x)) = \sum y^{\prime} \otimes h^{\prime}$, and $\Delta(k) = \sum k^{\prime} \otimes k^{\prime \prime}$.  Then
\begin{equation*}
\begin{split}
\Delta_K(m_s(x,k)) & = \bar{\Delta}(s(x)k) \\
  & \equiv s(x) \otimes k
\end{split}
\end{equation*}
modulo terms of the form $y^{\prime}k^{\prime} \otimes h^{\prime}k^{\prime \prime}$ with either $|y^{\prime}| < |s(x)|=n$, or $k^{\prime} \in I(K)$.  Otherwise said,
$$ \Delta_K(m_s(x,k)) \equiv s(x) \otimes k \mod (F_{n-1}M + I(K)M)\otimes I(H).$$
Thus
$$ \pi(\Delta_K(m_s(x,k))) \equiv x \otimes k \mod (Q_KM)^{<n}\otimes I(H),$$
and so both $\pi \circ \Delta_K \circ m_s$ and $\Delta_K \circ m_s$ are monic.
\end{proof}

\subsection{Proof of statements (a) and (b) of \thmref{prim thm}}

Let $B_r^{*,*} \subseteq Z_r^{*,*} \subseteq E_{r}^{*,*}$ denote the $r$--boundaries and $r$--cycles of the spectral sequence.

We can apply \lemref{handy lemma} to our spectral sequence by letting $H=H^*(C)$, $K=E_{\infty}^{0,*}$, and $M$ any of $E_r^{k,*}$, $Z_r^{k,*}$, $B_r^{k,*}$. We deduce

\begin{prop} \label{prim prop} For all $k$ and $r$, we have \\

\noindent{\bf (a)} $E_r^{k,*}$, $Z_r^{k,*}$, and $B_r^{k,*}$ are free $E_{\infty}^{0,*}$--modules, and \\

\noindent{\bf (b)} the composite $P_{H^*(C)}E_r^{k,*} \hra E_r^{k,*} \era Q_{E_r^{0,*}}E_r^{k,*}$ is monic.
\end{prop}

It follows that the short exact sequences of $E_{\infty}^{0,*}$--modules
$$ 0 \ra Z_r^{*,*} \ra E_r^{*,*} \ra B_r^{*,*} \ra 0, $$
and
$$ 0 \ra B_r^{*,*} \ra Z_r^{*,*} \ra E_{r+1}^{*,*} \ra 0, $$
are all split as $E_{\infty}^{0,*}$--modules.  Thus the spectral sequence remains a spectral sequence after applying $Q_{E_{\infty}^{0,*}}$.

Now let $A$ be a $(G,C)$--Duflot subalgebra of $H^*(G)$ as defined in \secref{main results section}: a subalgebra such that the composite $A \hra H^*(G) \xera{i^*} \im(i^*) = E_{\infty}^{0,*}$ is an isomorphism\footnote{We still need to show that such subalgebras exist.}.  We check the first two parts of \thmref{prim thm}:  (a) $H^*(G)$ is a free $A$--module so that the spectral sequence $\{Q_{E_{\infty}^{0,*}}E_r^{*,*}\}$ converges to $Q_AH^*(G)$, and (b) $P_{H^*(C)}H^*(G) \ra Q_AH^*(G)$ is monic.

Let $F_kBG$ be the inverse image of the $k$--skeleton of $BQ$ under the projection $BG \ra BQ$, and then let $F^k$ be the image of $H^*(BG) \ra H^*(F_kBG)$. Then $F^0 = E_{\infty}^{0,*}$, and for $k \geq 1$, there are short exact sequences
\begin{equation} \label{split sequence}
0 \ra E_{\infty}^{k,*} \ra F^k \ra F^{k-1} \ra 0
\end{equation}
of objects that are simultaneously $A$--modules and $H^*(C)$--comodules.  \propref{prim prop}(a) and induction on $k$ show these sequences split as $A$--modules.

Now consider the induced diagram
\begin{equation*}
\xymatrix{
0 \ar[r] & P_{H^*(C)}E_{\infty}^{k,*} \ar[d] \ar[r] & P_{H^*(C)}F^k \ar[d] \ar[r] & P_{H^*(C)}F^{k-1} \ar[d]  &  \\
0 \ar[r] & Q_AE_{\infty}^{k,*} \ar[r] & Q_AF^k \ar[r] & Q_AF^{k-1} \ar[r] & 0.  }
\end{equation*}
Here the top sequence is exact as indicated, as is the bottom, as (\ref{split sequence}) is split as $A$--modules.  The left vertical arrow is monic by \propref{prim prop}(b), as is the right vertical arrow, by induction on $k$, and it follows the middle arrow is also.

Thus we have proved that, for all $k$, $F^k$ is a free $A$--module and $P_{H^*(C)}F^k \ra Q_AF^k$ is monic.  As the connectivitity of the maps $H^*(G) \ra F^k$ goes to infinity as $k$ goes to infinity, we conclude that the same is true for $H^*(G)$.

\begin{rem}  The decreasing filtration of $H^*(G)$ induces a filtration on $P_CH^*(G)$ and $Q_AH^*(G)$.  We have shown that $Q_{E_{\infty}^{0,*}}E_{\infty}^{*,*}$ is the bigraded algebra associated to the filtration of $Q_AH^*(G)$.  By contrast, we can only conclude that the associated bigraded algebra of $P_CH^*(G)$ embeds in $P_{H^*(C)}E_{\infty}^{*,*}$.
\end{rem}

\subsection{Finite generation}  Statement (c) of \thmref{prim thm} says that both $P_CH^*(G)$ and $Q_AH^*(G)$ are finitely generated $H^*(Q)$--modules.  Our proof of this is similar to arguments used by L.Evens in \cite{evens}.

We first note that $i^*: H^*(G) \ra H^*(C)$ makes $H^*(C)$ into a finitely generated $H^*(G)$--module.  Otherwise put, $E_2^{0,*}$ is a finitely generated module over the ring $E_{\infty}^{0,*}$, which is Noetherian.

It follows that $E_2^{*,*} = H^*(Q) \otimes E_{2}^{0,*}$ is a finitely generated module over the Noetherian ring $H^*(Q) \otimes E_{\infty}^{0,*}$.  By induction on $r$, we conclude that, for all $r \geq 2$, $E_r^{*,*}$ is a finitely generated $H^*(Q) \otimes E_{\infty}^{0,*}$--module.

Passing to $E_{\infty}^{0,*}$--indecomposables, it follows that $Q_{E_{\infty}^{0,*}}E_{\infty}^{*,*}$ is a finitely generated $H^*(Q)$--module, and thus the same is true for $Q_AH^*(G)$, $P_{H^*(C)}E_{\infty}^{*,*}$, and $P_{H^*(C)}H^*(G)$.

\subsection{The image of inflation}

The quotient map $q: G \ra G/C$ induces the inflation homomorphism $q^*: H^*(G/C) \ra H^*(G)$. Its image, $\im(q^*)$, is an unstable subalgebra of $H^*(G)$ and also identifies with $E_{\infty}^{*,0}$ in the spectral sequence.

One approach to understanding $\im(q^*)$ is to try to understand $\ker(q^*)$.  Recall that the classifying homomorphism $\tau^{\#}: C^{\#} \ra H^2(G/C)$ corresponds to $d_2: E_2^{0,1} \ra E_2^{2,0}$.  Thus $\ker(q^*)$ is an ideal that is closed under Steenrod operations, and contains $\im(\tau^{\#})$.  As in \secref{main results section}, we let $I_{\tau} \subset H^*(G/C)$ be the smallest ideal with these properties.  Thus there is an epimorphism of unstable algebras $$H^*(G/C)/I_{\tau} \era \im(q^*),$$ which in many cases is an isomorphism.

As has already been said, $\im(q^*)$ is contained in the subalgebra $P_CH^*(G)$, but it seems worthwhile, at this point, to explicitly explain why.  The diagram
$$ C \times G
\begin{array}{c} m \\[-.08in] \longrightarrow \\[-.1in] \longrightarrow \\[-.1in] \pi
\end{array}
G \xra{q} G/C
$$
is a coequalizer diagram in the category of groups, i.e., a group homomorphism $f: G \ra H$ satisfies $f \circ m = f \circ \pi$ if and only if $f$ factors uniquely through $q$.  Applying cohomology, we have that $q^* \circ m^* = q^* \circ \pi^*$, and so
$$ \im(q^*) \subseteq \text{ Eq}\{m^*,\pi^*\} = P_CH^*(G).$$

In degree 1, inflation is as nice as possible.
\begin{lem} $q^*: H^1(G/C) \ra P_CH^1(G)$ is an isomorphism.
\end{lem}
\begin{proof}  The exact sequence arising from the corner of the spectral sequence,
$$ 0 \ra H^1(G/C) \xra{q^*} H^1(G) \xra{i^*} H^1(C),$$
can be viewed as the degree 1 part of a sequence of $H^*(C)$--comodules, if $H^*(G/C)$ is given the trivial comodule structure.  Taking primitives yields an exact sequence
$$ 0 \ra P_CH^1(G/C) \xra{q^*} P_CH^1(G) \xra{i^*} P_CH^1(C)$$
which identifies with
$$ 0 \ra H^1(G/C) \xra{q^*} P_CH^1(G) \ra 0,$$
as $P_CH^1(G/C)=H^1(G/C)$ and $P_CH^1(C) = 0$.
\end{proof}

In higher degrees, the inclusion $\im(q^*) \subseteq P_CH^*(G)$ certainly may be proper: see \exref{W(2) ex}.  However, we now show that the $\mathcal Nil_1$--closures of each of the maps
$$ H^*(G/C)/I_{\tau} \era \im(q^*) \hra P_CH^*(G)$$
is an isomorphism.  Equivalently, the composite is an $F$--isomorphism, as asserted in statement (d) of \thmref{prim thm}.

We prove this in two steps.

\begin{prop} $\displaystyle L_0(\im(q^*)) \simeq L_0(P_CH^*(G))$.
\end{prop}
\begin{proof}  By \propref{L0 prop}, we need to show that, for all $U$, there are bijections
$$ \Hom_{\K}(P_CH^*(G), H^*(U)) \simeq \Hom_{\K}(\im({\Inf}_{G/C}^G), H^*(U)).$$

Using that $H^*(U)$ is injective in $\K$, $\Hom_{\K}(P_CH^*(G), H^*(U))$ identifies with
$$ \text{Coeq } \{\Hom_{\K}(H^*(C\times G), H^*(U))
\begin{array}{c} m_* \\[-.08in] \longrightarrow \\[-.1in] \longrightarrow \\[-.1in] \pi_*
\end{array}
\Hom_{\K}(H^*(G), H^*(U))\},$$
and $\Hom_{\K}(\im(q^*), H^*(U))$ identifies with the image of
$$ \Hom_{\K}(H^*(G), H^*(U)) \xra{q^*} \Hom_{\K}(H^*(G/C), H^*(U)).$$

Lannes showed \cite[Prop.4.3.1]{L1} that
$\Rep(U,G) \simeq \Hom_{\K}(H^*(G), H^*(U))$,  where $\Rep(U,G)$ is set of orbits of $\Hom(U,G)$ under the  conjugation action of $G$.  Thus the next lemma is equivalent to the proposition.

\end{proof}
\begin{lem} The diagram of sets
$$ \Rep(U, C \times G)
\begin{array}{c} m_* \\[-.08in] \longrightarrow \\[-.1in] \longrightarrow \\[-.1in] \pi_*
\end{array}
\Rep(U,G) \xra{q_*} \Rep(U,G/C)$$
is exact in the following sense: given homomorphisms $\alpha, \beta: U \ra G$, $q \circ \alpha = q \circ \beta$ if and only if there exists $\gamma: U \ra C \times G$ such that $\alpha = m \circ \gamma$ and $\beta = \pi \circ \gamma$.
\end{lem}
\begin{proof} Such a $\gamma$ is equivalent to a pair $(\delta, \beta)$, where $\delta:U \ra C$ is a homomorphism satisfying $\delta(u)\beta(u)=\alpha(u)$ for all $u \in U$. Now suppose given $\alpha, \beta: U \ra G$ such that $q \circ \alpha = q \circ \beta$. Then the function $\delta: U \ra G$ defined by $\delta(u) = \alpha^{-1}(u)\beta(u)$ will take values in $C$ because $q \circ \alpha = q \circ \beta$, and will be a homomorphism because $C$ is central in $G$.
\end{proof}

\begin{prop} $\displaystyle L_0(H^*(G/C)/I_{\tau}) \simeq L_0(\im(q^*))$.
\end{prop}
\begin{proof}  In this case, $\Hom_{\K}(H^*(G/C)/I_{\tau}, H^*(U))$ identifies with
the set
$$ \{ \alpha \in \Rep(U,G/C) \ | \ \alpha^*(\tau) = 0 \in H^2(G/C;C) \},$$
while $\Hom_{\K}(\im(q^*), H^*(U))$ can be viewed as the set
$$ \{ \alpha \in \Rep(U,G/C) \ | \ \text{ $\alpha$ factors through $q:G \ra G/C$ }  \}.$$
But these sets are the same, because $\alpha^*(\tau)$ represents the top extension in the pullback diagram
\begin{equation*}
\xymatrix{
C \ar@{=}[d] \ar[r] & G(\alpha) \ar[d] \ar[r] & U \ar[d]^{\alpha} \\
C \ar[r] & G \ar[r]^q & G/C,}
\end{equation*}
and this extension is trivial if and only if $\alpha$ factors through $q$.
\end{proof}

We have a formula for $L_0(\im(q^*))$ analogous to the formula
$$L_0(H^*(G)) = \lim_{V \in \A(G)} H^*(V).$$  Let $\A(G,C)$ be the full subcategory of $\A(G)$ with objects the $V \in \A(G)$ containing $C$, and note that $H^*(V/C) = P_CH^*(V)$ for all such $V$.

\begin{prop} $\displaystyle L_0(\im(q^*)) \simeq \lim_{V \in \A(G,C)}H^*(V/C)$.
\end{prop}
\begin{proof} As $H^*(V/C)$ is $\mathcal Nil_1$--closed, so is $\displaystyle \lim_{V \in \A(G,C)}H^*(V/C)$.  Arguing as in the previous proofs, the proposition is equivalent to the statement that, for all $U$, the image of
$ \Rep(U,G) \xra{q_*} \Rep(U,G/C)$
identifies with $\displaystyle \colim_{V \in \A(G,C)}\Hom(U,V/C)$.  This is easily checked: details are left to the reader.
\end{proof}

These propositions allow us to quickly prove the last statement of \thmref{prim thm}: the Krull dimension of any of $H^*(G/C)/I_{\tau}$, $\im(q^*)$, $P_CH^*(G)$, or $Q_AH^*(G)$ equals (the $p$--rank of $G$) - (the rank of $C$).

To begin with, \propref{Krull dim prop} said that $\dim K = \dim L_0K$ if $K$ is a Noetherian unstable algebra. Thus $$ \dim H^*(G/C)/I_{\tau} = \dim \im(q^*)= \dim P_C H^*(G) = \dim \lim_{V \in \A(G,C)}H^*(V/C).$$
The second equality also follows from the fact that $P_CH^*(G)$ is a finitely generated $ H^*(G/C)/I_{\tau}$--module, and similarly $\dim H^*(G/C)/I_{\tau} = \dim Q_A H^*(G)$ is true.

As  $\displaystyle \prod _{V \in \A(G,C)}H^*(V/C)$ is a finitely generated $H^*(G/C)$--module, it certainly is also finitely generated over $\displaystyle \lim_{V \in \A(G,C)}H^*(V/C)$. Thus we conclude
\begin{equation*}
\begin{split}
\dim \lim_{V \in \A(G,C)}H^*(V/C) &
= \dim \prod _{V \in \A(G,C)}H^*(V/C) \\
  & = \max_{V \in \A(G,C)}\{\text{rank of }V/C\} \\
  & = \text{(the $p$--rank of $G$) - (the rank of $C$)}.
\end{split}
\end{equation*}

\section{Transgressions and the structure of $E_r^{0,*}$} \label{E0 section}

In this section we continue our examination of the spectral sequence associated to the central extension $C \xra{i} G \xra{q} Q$, where $C$ is $p$--elementary abelian of rank $c$.  We carefully describe the form of the differentials
$$ d_r: E_r^{0,r-1+*} \ra E_r^{r,*},$$
and prove that for all $r$, the Hopf algebra $E_r^{0,*} \subset H^*(C)$ must be a free commutative algebra\footnote{In the usual graded sense, if $p$ is odd.} of a standard form.  In particular, $\im(i^*) = E_{\infty}^{0,*}$ is free, and so a subalgebra of $H^*(G)$ generated by any lift of a miniminal set of generators of $\im(i^*)$ will be a $(G,C)$--Duflot subalgebra.

To begin our analysis, we know that $\tau^{\#}: C^{\#} \ra H^2(Q)$ corresponds to the transgression $ d_2: E_2^{0,1} \ra E_2^{2,0}$.

\subsection{The cokernel of $\tau$}

The cokernel of $\tau$ has a group theoretic meaning. The proof of the next lemma is left to the reader.

\begin{lem} Let $f: C \ra C_0$ be the cokernel of $\tau: H_2(Q) \ra C$.  Then $f$ factors through $C \xra{i} G$ and is the universal homomorphism from $C$ to a $p$--elementary abelian group with this property.
\end{lem}

As $f$ is split epic, the lemma tells us that $(G,C)$ is isomorphic to a pair of the form $C_0 \times (G_1, C_1)$, where no factor of $C_1$ splits off $G_1$. The spectral sequences will be related by
$ E_r^{*,*}(G,C) \simeq H^*(C_0) \otimes E_r^{*,*}(G_1, C_1)$, and $d_2: E_r^{0,1}(G_1, C_1) \ra E_r^{2,0}(G_1, C_1)$ will be monic.  Thus in our analysis, we can assume that $d_2: E_2^{0,1} \ra E_2^{2,0}$ is an inclusion if we wish.

\subsection{Properties of $d_r$}

Using that $E_r^{*,*}$ is free over $E_r^{0,*}$, we will see that the differentials we are interested in are `almost' determined by two standard properties.

The first is that $d_r$ is a derivation.

The second is the transgression theorem: recall that $E_r^{0,*}$ is an unstable subalgebra of $H^*(C)$ and $E_r^{*,0}$ is an unstable quotient algebra of $H^*(Q)$.  Given $x \in E_r^{0,r-1}$ and $a \in \A$,  $ax\in E_{r+|a|}^{0,r+|a|-1}$ and $d_{r+|a|}(ax)$ is represented by  $ad_r(x)$.

\subsection{The kernel of the DeRham and Koszul derivations}

Our differentials are modelled by two standard derivations.

Let $S^*(V)$ be the symmetric algebra on a graded $\F_p$--vector space $V$, with $V$ concentrated in even degrees if $p$ is odd.  The DeRham derivation is the derivation
$$d_V: S^*(V) \ra S^*(V) \otimes \Sigma V$$
determined by letting $d_V(v) = 1 \otimes \sigma v$ for $v \in V$.

We will need to know its kernel.  Let $\Phi: S^*(V) \ra S^*(V)$ denote the $p^{th}$ power map.

\begin{lem} \label{DeRham lem} The kernel of $d_V$ is $S^*(\Phi(V))$.
\end{lem}

This is a special case of a result due to Cartier \cite{cartier}.  For completeness, we sketch an elegant argument we learned from \cite[proof of Prop.3.3]{fls}.  The kernel of $d_V$ is $H^0$ of the DeRham complex $\Omega^*(V) = (S^*(V) \otimes \Lambda^*(\Sigma V), d_V)$.  Since $\Omega^*(V\oplus W) \simeq \Omega^*(V) \otimes \Omega^*(W)$, the Kunneth theorem allows one to reduce to the case when $V$ is one dimensional, where is it easily checked.

Let $\Lambda^*(V)$ be the exterior algebra on a graded $\F_p$--vector space $V$, with $V$ concentrated in odd degrees. The Koszul derivation is the derivation
$$\delta_V: \Lambda^*(V) \ra \Lambda^*(V) \otimes \Sigma V$$
determined by letting $\delta_V(v) = 1 \otimes \sigma v$ for $v \in V$.  This is the bottom of the Koszul complex $(\Lambda^*(V) \otimes S^*(\Sigma V), \delta_V)$, which is acyclic, and we have

\begin{lem} \label{Koszul lem} The kernel of $\delta_V$ is $\F_p$, i.e. $\delta_V$ is monic in positive degrees.
\end{lem}

\subsection{The structure of $E_r^{0,*}$ when $p=2$}

If $p=2$, $E_2^{0,*} = S^*(C^{\#})$.  As the squaring operation $\Phi$ in degree $n$ corresponds to the Steenrod operation $Sq^n$,  the image of $\Phi^{k-1}: E_2^{0,1} \ra E_2^{0,2^k}$ lands in the subspace $E_{2^k+1}^{0,2^k}$.  Thus we can define an increasing filtration of $C^{\#}$,
$$ C_0^{\#} \subseteq C_1^{\#} \subseteq C_2^{\#} \subseteq \dots$$
by letting $C_k^{\#}$ be the kernel of the composite
$$ E_2^{0,1} \xra{\Phi^{k}} E_{2^k+1}^{0,2^k} \xra{d_{2^k+1}} E_{2^k+1}^{2^k+1,0}.$$

\begin{thm} \label{E0 structure thm} The only possible nonzero differentials
$$ d_r: E_r^{0,*} \ra E_r^{r,*}$$
are $d_{2^k+1}$ with $k=0,1,2,\dots$. For each such $k$, $ E_{2^k+2}^{0,*} = \ker d_{2^k+1}$ equals
$$ S^*(C_0^{\#} + \Phi(C_1^{\#}) + \dots + \Phi^{k}(C_{k}^{\#})+ \Phi^{k+1}(C^{\#})).$$
\end{thm}
This polynomial algebra is noncanonically isomorphic to
$$ S^*(V_0 \oplus \Phi(V_1) \oplus \dots \oplus \Phi^{k}(V_{k}) \oplus \Phi^{k+1}(C^{\#}/C_{k}^{\#})),$$
where $V_i = C_{i}^{\#}/C_{i-1}^{\#}$.  In particular, one can choose a basis for $C^{\#}$ such that $E_{\infty}^{0,*}$ has the form described in \secref{res subsection}.

\begin{rem}  In practical terms, the filtration of $C^{\#}$ is often determined by the extension homomorphism $\tau: H_2(Q) \ra C$ together with the action of the Steenrod operations on $H^*(Q)$.

Note that $C_k^{\#}$ will also be the kernel of the composite
$$ E_2^{0,1} \xra{d_2} E_{2}^{2,0} \xra{Sq^1}  \dots \xra{Sq^{2^{k-1}}} E_{2}^{2^k+1,0} \era  E_{2^k+1}^{2^k+1,0},$$
where the last map is the evident quotient. Equivalently, $C_k^{\#}$ is the kernel of
$$ C^{\#} \xra{\tau_G^{\#}} H^2(Q) \xra{Sq^1}  \dots \xra{Sq^{2^{k-1}}} H^{2^k+1}(Q) \era E_{2^k+1}^{2^k+1,0}.$$
Let $I_{\tau}(k) \subset H^*(Q)$ be the ideal generated by $\A(k-1)\cdot \im(\tau^{\#})$, where $\A(k) \subset \A$ is the subalgebra generated by $Sq^1, \dots, Sq^{2^k}$.  Then the quotient map $H^{*}(Q) \era E_{2^k+1}^{*,0}$ factors
$$ H^{*}(Q) \era H^{*}(Q)/I_{\tau}(k-1) \era E_{2^k+1}^{*,0},$$
and the second map is often an isomorphism.

\end{rem}

\begin{proof}[Proof of \thmref{E0 structure thm}]

By definition, $\Phi^k(C_k^{\#})$ is the kernel of the transgression $d_{2^k+1}: \Phi^k(C^{\#}) \ra E_{2^k+1}^{2^k+1,0}$, and so consists of permanent cycles. It follows that, for each $k$, the subalgebra $S(k) = S^*(C_0^{\#} + \Phi(C_1^{\#}) + \dots + \Phi^{k}(C_{k}^{\#}))$ is also contained in $E_{\infty}^{0,*}$.

By induction, we now show that
$$ E_{2^k+1}^{0,*} = S^*(C_0^{\#} + \Phi(C_1^{\#}) + \dots + \Phi^{k-1}(C_{k-1}^{\#})+ \Phi^{k}(C^{\#})).$$
When $k=0$, this just says that $E_2^{0,*} = S^*(C^{\#})$, and so is certainly true.  Now we assume the statement for $k$ and prove it with $k$ replaced by $k+1$.

Let $V = \Phi^k(C^{\#})/\Phi^k(C_k^{\#})$.   Then $E_{2^k+1}^{0,*} \simeq S(k) \otimes S^*(V)$, $\Sigma V$ identifies with the image of $d_{2^k+1}: \Phi^k(C^{\#}) \ra E_{2^k+1}^{2^k+1,0}$, and we have a commutative diagram
\begin{equation*}
\xymatrix{
E_{2^k+1}^{0,*} \ar[d]^{\wr} \ar[r]^{d_{2^k+1}} & E_{2^k+1}^{2^k+1,*}   \\
S(k) \otimes S^*(V) \ar[r]^-{1 \otimes d_V} & S(k) \otimes S^*(V) \otimes \Sigma V.  \ar[u]}
\end{equation*}
Here the right vertical map is induced from the inclusion $\Sigma V \subseteq E_{2^k+1}^{2^k+1,0}$ using the $E_{2^k+1}^{0,*}$--module structure on $E_{2^k+1}^{2^k+1,*}$.

We have reached the key point in our proof: as in our proof of \thmref{prim thm}, \lemref{handy lemma} shows that this module structure is free, and thus the right vertical map is monic.  It follows that the kernel of the top map identifies with $\ker(1 \otimes d_V)$ which equals $S(k) \otimes S^*(\Phi(V))$, by \lemref{DeRham lem}.  Otherwise said,
$$ E_{2^k+2}^{0,*} = S^*(C_0^{\#} + \Phi(C_1^{\#}) + \dots + \Phi^{k}(C_{k}^{\#})+ \Phi^{k+1}(C^{\#})).$$

Finally, we note that once we know that $ E_{2^k+2}^{0,*}$ has this form, $E_{2^{k+1}+1}^{0,*} = E_{2^k+2}^{0,*}$ follows by the transgression theorem.
\end{proof}

\subsection{The structure of $E_r^{0,*}$ when $p$ is odd}

When $p$ is odd, $$E_2^{0,*} = \Lambda^*(C^{\#}) \otimes S^*(\beta(C^{\#})).$$  As the $p^{th}$ power operation $\Phi$ in degree $2n$ corresponds to the Steenrod operation $\mathcal P^n$, the image of $\Phi^{k}\circ \beta: E_2^{0,1} \ra E_2^{0,2p^k}$ lands in the subspace $E_{2p^k+1}^{0,2p^k}$.  Thus we can define an increasing filtration of $C^{\#}$,
$$ C_0^{\#} \subseteq C_1^{\#} \subseteq C_2^{\#} \subseteq \dots$$
by letting $C_0^{\#}$ be the kernel of $E_2^{0,1} \xra{d_2} E_2^{2,0}$, and, for $k \geq 0$,
$C_{k+1}^{\#}$ be the kernel of the composite
$$ E_2^{0,1} \xra{\Phi^{k}\circ \beta} E_{2p^k+1}^{0,2p^k} \xra{d_{2p^k+1}} E_{2p^k+1}^{2p^k+1,0}.$$

\begin{thm} \label{E0 structure thm, p odd} The only possible nonzero differentials
$$ d_r: E_r^{0,*} \ra E_r^{r,*}$$
are $d_2$ and $d_{2p^k+1}$ with $k=0,1,2,\dots$.  Furthermore,
$$ E_3^{0,*} = \Lambda^*(C_0^{\#})\otimes  S^*(\beta(C^{\#})),$$
and, for each $k \geq 0$,
$$ E_{2p^k+2}^{0,*} = \Lambda^*(C_0^{\#})\otimes  S^*(\beta(C_1^{\#}) + \Phi\beta(C_2^{\#}) + \dots + \Phi^{k}\beta(C_{k+1}^{\#})+ \Phi^{k+1}\beta(C^{\#})).$$
\end{thm}
This free commutative algebra is noncanonically isomorphic to
$$ \Lambda^*(V_0)\otimes S^*(\beta(V_0) \oplus \beta(V_1) \oplus \dots \oplus \Phi^{k}\beta(V_{k+1}) \oplus \Phi^{k+1}\beta(C^{\#}/C_{k+1}^{\#})),$$
where $V_i = C_{i}^{\#}/C_{i-1}^{\#}$.  In particular, one can choose a basis for $C^{\#}$ such that $E_{\infty}^{0,*}$ has the form described in \secref{res subsection}.

\begin{rem}  Similar to the case when $p=2$, $C_{k+1}^{\#}$ will be the kernel of
$$ C^{\#} \xra{\tau^{\#}} H^2(Q) \xra{\beta} H^3(Q) \xra{\mathcal P^{1}} \dots \xra{\mathcal P^{p^{k-1}}} H^{2p^k+1}(Q) \era E_{2p^k+1}^{*,0},$$
and the quotient map $H^{*}(Q) \era E_{2p^k+1}^{*,0}$ factors as $$ H^{*}(Q) \era H^{*}(Q)/I_{\tau}(k) \era E_{2p^k+1}^{*,0}$$
where $I_{\tau}(k) \subset H^*(Q)$ is the ideal generated by $\A(k-1)\cdot \im(\tau^{\#})$.  Here $\A(k) \subset \A$ is the subalgebra generated by $\beta, \mathcal P^1, \mathcal P^p, \dots, \mathcal P^{p^{k-1}}$.
\end{rem}

\begin{proof}[Proof of \thmref{E0 structure thm, p odd}]  We just sketch the proof, as it follows along the lines of the proof of the $p=2$ version of the theorem.

To compute $E_3^{0,*}$, let $V = C^{\#}/C_0^{\#}$. Then $\ker d_2$ identifies with the kernel of
$$ \Lambda^*(C_0^{\#}) \otimes S^*(\beta (C^{\#})) \otimes \Lambda^*(V) \xra{1 \otimes 1 \otimes \delta_V} \Lambda^*(C_0^{\#}) \otimes S^*(\beta (C^{\#})) \otimes \Lambda^*(V) \otimes \Sigma V.$$
The formula for $E_3^{0,*}$ thus follows from \lemref{Koszul lem}.

To compute $E_{2p^k+2}^{0,*}$ for $k \geq 0$, let $V = \Phi^{k}\beta (C^{\#})/\Phi^{k}\beta(C_{k+1}^{\#})$.  Then the subalgebra
$$S(k) = \Lambda^*(C_0^{\#})\otimes  S^*(\beta(C_1^{\#}) + \Phi\beta(C_2^{\#}) + \dots + \Phi^{k}\beta(C_{k+1}^{\#}))$$
is all permanent cycles.  Using that  $E_{2p^k+1}^{2p^k+1,*}$ is a free $E_{2p^k+1}^{0,*}$--module, $\ker d_{2p^k+1}$ identifies with the kernel of
$$ S(k) \otimes S^*(V) \xra{1 \otimes d_V} S(k) \otimes S^*(V) \otimes \Sigma V,$$
and the formula for $E_{2p^k+2}^{0,*}$ follows from \lemref{DeRham lem}.
\end{proof}

\section{$p$--central groups} \label{pC section}  In this section we prove our main theorems about $p$--central groups: \thmref{PC top prim thm} and \thmref{PC d0 thm}.

We begin by recalling some notation from \secref{main results section}.  If $C=C(G)$ is the maximal central $p$--elementary abelian subgroup of a finite group $G$, we have shown that $C^{\#}=H^1(C)$ admits an ordered basis $(x_1, \dots, x_c)$ so that,  if $y_j = \beta (x_j)$ for $p$ odd,
\begin{equation*}
\Res_C^G(H^*(G)) =
\begin{cases}
\F_2[x_1^{2^{j_1}},\dots,x_c^{2^{j_c}}] & \text{if $p=2$}  \\

\F_p[y_1^{p^{j_1}},\dots,y_b^{p^{j_b}},y_{b+1},\dots,y_c] \otimes \Lambda(x_{b+1},\dots,x_c)  & \text{if $p$ is odd},
\end{cases}
\end{equation*}
with the $j_i$ forming a sequence of nondecreasing nonnegative integers.

Then we say that $G$ has type $[a_1,\dots,a_c]$ where
\begin{equation*}
(a_1,\dots,a_c) =
\begin{cases}
(2^{j_1}, \dots, 2^{j_c}) & \text{if $p=2$}  \\ (2p^{j_1}, \dots, 2p^{j_b},1, \dots,1) & \text{if $p$ is odd},
\end{cases}
\end{equation*}
and we let $\displaystyle e(G) = \sum_{i=1}^c (a_i-1)$ and $h(G) =
\begin{cases}
2p^{k-1} & \text{if } a_1 = 2p^k \text{ with } k \geq 1 \\
1 & \text{if } a_1 = 2 \\ 0 & \text{if } a_1=1.
\end{cases}$

We have the following lemma about products.

\begin{lem}  Suppose $G_0$ and $G_1$ have maximal central $p$--elementary abelian subgroups $C_0$ and $C_1$, and Duflot subalgebras $A_0$ and $A_1$.  Then the following hold.  \\

\noindent{\bf (a)} $C_0 \times C_1 = C(G_0 \times G_1)$, and $$P_{C_0 \times C_1}H^*(G_0 \times G_1) = P_{C_0}H^*(G_1) \otimes P_{C_1}H^*(G_1).$$

\noindent{\bf (b)} $A_0 \otimes A_1$ will be a Duflot subalgebra for $G_0 \times G_1$, and
$$Q_{A_0 \otimes A_1}H^*(G_0 \times G_1) = Q_{A_0}H^*(G_1) \otimes Q_{A_1}H^*(G_1).$$

\noindent{\bf (c)} $e(G_0 \times G_1) = e(G_0) + e(G_1)$, and $h(G_0 \times G_1)= \max\{h(G_0), h(G_1)\}$.
\end{lem}

Note that a subgroup $H$ of a $p$--central group $G$ is again $p$--central, and $C(H) = C(G) \cap H$.  The next lemma is easily deduced.

\begin{lem} Let $G$ be $p$--central, and let $A \subset H^*(G)$ be a Duflot subalgebra. If $j:H < G$ is a subgroup, then $e(H) \leq e(G)$, $h(H) \leq h(G)$, and $j^*(A)$ will be a Duflot subalgebra of $H^*(H)$.
\end{lem}

Thanks to this lemma, \corref{PC d0 cor} immediately follows from \thmref{PC d0 thm}.

\begin{rem}  The example $H=\Z/4 < \Z/8 = G$ shows that the inequalities of the lemma can be equalities, even when $H$ is a proper subgroup of a $p$--group $G$.
\end{rem}

\subsection{Benson--Carlson duality}   If $G$ is $p$--central, then  $Q_AH^*(G)$ will be a finite dimensional $\F_p$--algebra if $A$ is any Duflot subalgebra.  Benson and Carlson tell us much more:

\begin{thm} \label{duality thm}  If $G$ is $p$--central and $A$ is a Duflot subalgebra of $H^*(G)$, then $Q_AH^*(G)$ is a Poincar\'e duality algebra with top class in degree $e(G)$.
\end{thm}

Under the assumption that $A$ is a polynomial algebra (always true if $p=2$), this is an immediate application of the main theorem in \cite{bc}.  The general case reduces to this one: $G$ and $A$ will admit decompositions $G = C_0 \times G_1$ and $A = H^*(C_0) \otimes A_1$, with $C_0$ $p$--elementary, $G_1$ having no $\Z/p$ summands, and $A_1$ a (necessarily polynomial) Duflot subalgebra of $H^*(G_1)$.  Then $Q_AH^*(G) = Q_{A_1}H^*(G_1)$, and $e(G) = e(G_1)$.

\subsection{Proof of \thmref{PC top prim thm}}  Let $G$ be $p$--central, $C=C(G)$, and $A \subset H^*(G)$ a Duflot subalgebra.  We now prove the various parts of \thmref{PC top prim thm}.

Firstly, \thmref{duality thm} implies that $Q_AH^*(G)$ is zero in degrees greater than $e(G)$, and one dimensional in degree $e(G)$.

Now consider the Serre spectral sequence for $C \ra G \ra G/C$, as studied in \thmref{prim thm}.  The bigraded algebra $Q_{E_{\infty}^{0,*}}E_{\infty}^{*,*}$ is the graded object associated to a decreasing filtration of the Poincar\'e duality algebra $Q_AH^*(G)$ with top degree $e(G)$.  This forces the following to be true: there is a largest $s$, $s(G)$, such that $E_{\infty}^{s,*}$ is nonzero, $Q_{E_{\infty}^{0,*}}E_{\infty}^{s(G),*}$ will be one dimensional and concentrated in total degree $e(G)$, and nonzero classes in $E_{\infty}^{s(G),e(G)-s(G)} \subset H^{e(G)}(G)$ will be Poincar\'e duality classes.

These classes will also be $H^*(C)$--comodule primitives, as  $E_{\infty}^{s(G),*}$ is a sub--$H^*(C)$--comodule of $H^*(G)$, and everything in lowest degree must be primitive.  As $P_CH^*(G)$ is contained in $Q_AH^*(G)$, we conclude that $P_CH^*(G)$ is also zero in degrees greater than $e(G)$, and one dimensional in degree $e(G)$.

By \corref{LF cor}, $P_CH^{e(G)}H^*(G)$ is also be the top nonzero degree of $H^*(G)_{LF}$, and so consists of classes annihilated by all positive degree Steenrod operations.

It remains to show that, under the additional assumption that $G$ is a $p$--group, $P_CH^{e(G)}H^*(G)$ is essential cohomology.  This we prove in the next subsection.

\subsection{$p$--central $p$--groups and essential cohomology}

Let $P$ be a $p$--central $p$--group.  We have shown that $H^{e(P)}(P)_{LF} = P_{C(P)}H^{e(P)}(P)$ is a one dimensional subspace of $H^{e(P)}(P)$.

\begin{prop} $H^{e(P)}(P)_{LF}$ is essential.
\end{prop}
\begin{proof} As $P$ is a $p$--group, maximal proper subgroups have the form $j:Q<P$, where $Q$ is the kernel of a nonzero homomorphism $x: P \ra \Z/p$.  We need to show that $j^*(\zeta) = 0 \in H^*(Q)$ if $\zeta \in H^{e(P)}(P)_{LF}$ is nonzero.

The map $j^*: H^*(P) \ra H^*(Q)$ will take $H^{e(P)}(P)_{LF}$ to $H^{e(P)}(Q)_{LF}$.  If $e(Q)<e(P)$, we are done: $j^*(\zeta)$ will be an element of a zero group.

If $e(Q)=e(P)$, we reason as follows.  Let $A$ be a Duflot subalgebra of $H^*(G)$, so that $j^*(A)$ is a Duflot subalgebra of $H^*(Q)$.  If $j^*(\zeta) \neq 0$, it will project to a nonzero element in $Q_{j^*(A)}H^*(Q)$. We show that this is impossible.  Regard $x$ as an nonzero element in $H^1(P)$.  By construction, $j^*(x) = 0 \in H^1(Q)$.  By Poincar\'e duality, there exists $y \in H^*(P)$ such that $\zeta = xy \in Q_AH^*(P)$. But then $j^*(\zeta) = j^*(x)j^*(y) = 0 \in Q_{j^*(A)}H^*(Q)$.
\end{proof}

Let $A(P,P)$ be the two sided Burnside ring over $\F_p$: the $\F_p$--algebra with basis given by equivalence classes of diagrams $P \geq Q \xra{\alpha} P$, and multiplication defined using the double coset formula\footnote{There are more elegant descriptions, but this is better for our purposes.}.  If $J$ is the ideal generated by all such diagrams with $\alpha$ {\em not} an isomorphism, then $A(P,P)/J \simeq \F_p[Out(P)]$, the group ring of the outer automorphism group.

Using transfers (a.k.a. induction), $A(P,P)$ acts on $H^*(P)$, with a basis element $[P \geq Q \xra{\alpha} P]$ inducing $H^*(P) \xra{\alpha^*} H^*(Q) \xra{Tr_Q^P} H^*(P)$. As these are unstable $\A$--module maps, it follows that $H^{e(P)}(P)_{LF}$ is a one dimensional $A(P,P)$--submodule.

\begin{cor} \label{J cor} The ideal $J$ acts trivially on $H^{e(P)}(P)_{LF}$.
\end{cor}
\begin{proof}  The previous proposition shows that if a homomorphism $\alpha: Q \ra P$ is not onto, then $\alpha^*(H^{e(P)}(P)_{LF}) = 0$.
\end{proof}

It follows that the $A(P,P)$--module $H^{e(P)}(P)_{LF}$ is the pullback of a one dimensional representation of $Out(P)$ over the prime field $\F_p$.  We let $\omega(P)$ denote this representation.  Clearly $\omega(P)$ will be trivial if $p=2$, but this need not be the case when $p$ is odd.

\begin{ex}  Let $p=3$.  Then $\omega(\Z/9) = H^1(\Z/9)$ is nontrivial, as $-1: \Z/9 \ra \Z/9$ induces multiplication by $-1$ on $H^1(\Z/9)$.
\end{ex}

\subsection{ $d_0(G)$ when $G$ has a $p$--central $p$--Sylow subgroup}  \label{d0(P) subsection} We prove the parts of \thmref{PC d0 thm} involving $d_0$.

Firstly, if $G$ is $p$--central, then \corref{Rd cor} says that
$$\displaystyle \bar R_dH^*(G) \simeq  H^*(C(G)) \otimes P_{C(G)}H^d(G).$$  Since $d_0(G)$ is the largest $d$ such that $\bar R_dH^*(G) \neq 0$, it follows that $d_0(G)$ will equal the top nonzero degree of $P_{C(G)}H^*(G)$, which we have computed to be $e(G)$.

Now suppose that $G$ is not necessarily $p$--central, but has a $p$--central $p$--Sylow subgroup $P$. We show that then $d_0(G) = d_0(P)$.

We need to show that the largest $d$ such that $\bar R_dH^*(G) \neq 0$ is $d=d_0(P)=e(P)$. Let $e_1 \in A(P,P)$ be an idempotent chosen so that $A(P,P)e_1$ is the projective cover of $\epsilon$, the trivial $\F_p[Out(P)]$--module, pulled back to $A(P,P)$.  Standard arguments show that there are inclusions
$$ e_1 \bar R_dH^*(P) \subseteq \bar R_dH^*(G) \subseteq \bar R_dH^*(P).$$
Thus it suffices to show that $ e_1 \bar R_{e(P)}H^*(P) \neq 0$. Otherwise said, it suffices to show that $\epsilon$ is a composition factor in the $A(P,P)$--modules $\bar R_{e(P)}H^*(P)$.

If $p=2$, we are done: by \corref{J cor}, $\bar R_{e(P)}H^{e(P)}(P) = H^{e(P)}(P)_{LF} \simeq \epsilon$.  As a bonus, we learn that $H^*(G)_{LF}$ is one dimensional in degree $e(P)$.

When $p$ is odd, more care (and maybe luck) is needed.  Recall that $\bar R_{e(P)}H^*(P) = H^*(C(P)) \otimes H^{e(P)}(P)_{LF}$.  The fact that $J$ acts as 0 on $H^{e(P)}(P)_{LF}$ implies that same is true for $H^*(C(P)) \otimes H^{e(P)}(P)_{LF}$.  Thus we just need to show that the trivial $Out(P)$--module occurs as a composition factor in $H^*(C(P)) \otimes H^{e(P)}(P)_{LF} = H^*(C(P)) \otimes \omega(P)$, or, equivalently, that $\omega(P)^{-1}$ occurs as an $Out(P)$-composition factor $H^*(C(P))$.  We are done with the following lemma\footnote{This lemma is false if $p=2$, as the example $P=Q_8$ illustrates.}.

\begin{lem} If $P$ is a $p$--central $p$--group with $p$ an odd prime, then every irreducible $\F_p[Out(P)]$--module occurs as a composition factor of $H^*(C(P))$.
\end{lem}

The lemma, in a stronger form than stated, follows by combining  \cite[Prop.5.7 and Cor.6.8]{k4}.  The key point is that, since $C(P) = \Omega_1(P)$, the kernel of $Aut(P) \ra Aut(C(P))$ will be a $p$--group if $p$ is odd \cite[Thm.5.3.10]{gor}.

\begin{ex}  Let $p=3$ and $G$ be the semidirect product $\Z/9 \rtimes \Z/2$.  Then $d_0(G) = d_0(\Z/9) = 1$, but $\tilde{H}^*(G)_{LF} = \mathbf 0$.
\end{ex}

\subsection{$d_1(G)$ when $G$ is $p$--central}  In this subsection, let $G$ be $p$--central.  We show that $d_1(G) = e(G) + h(G)$.

We get control of $d_1(G)$ by working directly with the desuspended composition factors $R_dH^*(G)$ of $H^*(G)$, rather than their $\mathcal Nil_1$--localizations $\bar R_dH^*(G)$, as was done in our calculation of $d_0(G)$.

To simplify notation, write $R_d$ for $R_dH^*(G)$, $\bar R_d$ for $\bar R_dH^*(G)$, and $C$ for $C(G)$.  We have that $\bar R_0 = H^*(C)$, and $R_0 = \im(i^*)$, where $i: C \hra G$ is the inclusion.

In the nilpotent filtration of $H^*(G)$, the last nonzero submodule, $nil_{e(G)}H^*(G) = \Sigma^{e(G)}R_{e(G)}$, has been shown to be isomorphic to $\Sigma^{e(G)}R_0$ as an unstable module.  Thus $d_1$ of this submodule of $H^*(G)$ equals $e(G)+d_1(R_0)$.  By \lemref{d1 prop 2}, the next lemma implies that if $d<e(G)$, $d_1(\Sigma^dR_d)$ is strictly smaller than this.

\begin{lem} \label{decomp lemma}  Each $R_d$ with $d < e(G)$ admits a filtration by unstable modules with subquotients all of the form $\Sigma^kR_0$ with $d+k<e(G)$.
\end{lem}

Again appealing to \lemref{d1 prop 2}, we then have the next corollary.

\begin{cor} $d_1(G) = d_1(\Sigma^{e(G)}R_{e(G)}) = e(G) + d_1(R_0)$.
\end{cor}

Thus we will have proved that, when $G$ is $p$--central, $d_1(G) = e(G) + h(G)$, once we have proved \lemref{decomp lemma}, and calculated that $d_1(R_0) = h(G)$.  We begin the proof of \lemref{decomp lemma} here, and then both finish it, and calculate $d_1(R_0)$, in the two subsections that follow, which correspond to the cases $p=2$ and $p$ odd.

\begin{proof}[Proof of \lemref{decomp lemma}] For all $d$, we have inclusions
$$ R_0 \otimes P_CH^d(G) \subseteq R_d \subseteq \bar R_0 \otimes P_CH^d(G),$$
where $P_CH^d(G)$ is regarded as an unstable module concentrated in degree $0$.  These are inclusions of unstable modules, enriched with compatible $R_0$--module structures and $\bar R_0$--comodule structures.  Call the category of such objects $R_0-\bar R_0-\U$.

Say that $M \in R_0-\bar R_0-\U$ admits a nice filtration if it admits a filtration in $R_0-\bar R_0-\U$ with subquotients all of the form $\Sigma^kR_0$.  We will show that each $R_d$ admits a nice filtration.

(That the composition factors will then also satisfy $d+k<e(G)$ follows immediately from the fact that $Q_{R_0}(\Sigma^d R_*)$ is a graded object associated to $Q_AH^*(G)$, which we know is one dimensional in degree $e(G)$ and zero above that.)

We claim that, if $N$ admits a nice filtration, and $M \subseteq N$, then $M$ also admits a nice filtration.  To see this, suppose $F_0N \subseteq F_1N \subseteq \dots$ is a filtration of $N$ with $F_jN/F_{j-1}N = \Sigma^{k_j}R_0$.  Let $F_jM = M \cap F_jN$. Then $F_jM/F_{j-1}M \subseteq F_jN/F_{j-1}N = \Sigma^{k_j}R_0$ will be an inclusion of objects in $R_0-\bar R_0-\U$ that will be split as $R_0$--modules, thanks to \corref{K-H-cor}.  We conclude that $F_jM/F_{j-1}M$ is either 0 or $\Sigma^{k_j}R_0$.

Thus to prove $R_d$ has a nice filtration, it suffices to prove that $\bar R_0 \otimes P_CH^d(G)$ has a nice filtration, or just that $\bar R_0$ has a nice filtration.  We  show this in the next two subsections, which separately deal with the cases $p=2$ and $p$ is odd.
\end{proof}

\subsection{A calculation of $d_1(R_0)$, and a nice filtration of $\bar R_0$, when $p=2$}

Suppose that $p=2$, and that
$$R_0 = \F_2[x_1^{2^{j_1}}, \dots, x_1^{2^{j_c}}] \subseteq \F_2[x_1, \dots, x_c] = \bar R_0,$$
with $j_1 \geq  \dots \geq j_c$.
We show the following.

\begin{lem} \label{good filt lemma, p=2} $\bar R_0$ has a good filtration as an object in $R_0-\bar R_0 - \U$.
\end{lem}

\begin{lem} \label{d_1R_0, p=2} If $j_1>0$, $d_1(R_0) = 2^{{j_1}-1}$.
\end{lem}

\begin{proof}[Proof of \lemref{good filt lemma, p=2}]
For all $1\leq b \leq c$ and $1 \leq i_b \leq j_b$, the module
$$R(i_1, \dots, i_c) = \F_2[x_1^{2^{i_1}}, \dots, x_c^{2^{i_c}}]$$ will be an object in $R_0-\bar R_0-\U$ in the evident way.

Clearly $R(j_1, \dots, j_c)=R_0$ admits a nice filtration.  The short exact sequences
$$ 0 \ra R(i_1, \dots, i_c) \ra R(i_1, \dots,i_b-1, \dots, i_c) \ra \Sigma^{2^{i_b}}R(i_1, \dots, i_c) \ra 0$$
then shows that if $R(i_1, \dots, i_c)$ admits a nice filtration, so does $R(i_1, \dots,i_b-1, \dots, i_c)$.  By downward induction, we conclude that $R(0,\dots,0) = \bar R_0$ admits a nice filtration.
\end{proof}

\begin{proof}[Proof of \lemref{d_1R_0, p=2}]  By \propref{d1 prop 1}(c), it suffices to prove that, if $j>0$,
$$ d_1(\F_2[x^{2^j}]) = 2^{j-1}.$$

In the short exact sequence
$$ 0 \ra \F_2[x^{2^j}] \ra \F_2[x^{2^{j-1}}] \ra \Sigma^{2^{j-1}}\F_2[x^{2^j}] \ra 0,$$
$d_1$ of the middle term is strictly less than $d_0(\Sigma^{2^{j-1}}\F_2[x^{2^j}]) = 2^{j-1}$: this is clear when $j=1$, and for larger $j$ this follows by an inductive hypothesis.  Thus \propref{d1 prop 2}(b) applies to say that $ d_1(\F_2[x^{2^j}]) = 2^{j-1}.$
\end{proof}

\subsection{A calculation of $d_1(R_0)$, and a nice filtration of $\bar R_0$, when $p$ is odd}

Suppose that $p$ is odd.  We can assume that
$$R_0 = \F_p[y_1^{p^{j_1}}, \dots, y_1^{p^{j_c}}] \subseteq \Lambda^*(x_1, \dots, x_c) \otimes \F_p[y_1, \dots, y_c] = \bar R_0,$$
with $j_1 \geq  \dots \geq j_c$, and we show the following.

\begin{lem} \label{good filt lemma, p odd} $\bar R_0$ has a good filtration as an object in $R_0-\bar R_0 - \U$.
\end{lem}

\begin{lem} \label{d_1R_0, p odd} If $j_1=0$, $d_1(R_0) = 1$.  If $j_1>0$, $d_1(R_0) = 2p^{{j_1}-1}$.
\end{lem}

\begin{proof}[Proof of \lemref{good filt lemma, p odd}]
As a first step, we note that the filtration of $\bar R_0$ given by letting $F_k\bar R_0 = \Lambda^{\leq k}(x_1, \dots, x_c) \otimes \F_p[y_1, \dots, y_c]$ is a filtration in the category $R_0-\bar R_0 - \U$, and the associated subquotients are direct sums of suspensions of $\F_p[y_1, \dots, y_c]$.  It follows that it suffices to prove the lemma with $\bar R_0$ replaced by $\F_p[y_1, \dots, y_c]$.

Our next reduction will allow us to reduce to the case when $c=1$.

If $K$ is a subHopf algebra of a Hopf algebra $H$, and both objects and all structure maps are in $\U$, one has a category $K-H-\U$, analogous to $R_0-\bar R_0-\U$.  One can then say that $M \in K-H-\U$ has a good filtration if it has a filtration with subquotients that are all suspensions of $K$.  It is easy to see that if $M_1 \in K_1-H_1-\U$ and $M_2 \in K_2-H_2-\U$ has a good filtration, then so does $M_1\otimes M_2$, viewed as an object in $K_1\otimes K_2-H_1\otimes H_2 -\U$.

Applying this observation to the evident tensor decompositions of $K=R_0$ and $H=\F_p[y_1, \dots, y_c]$, we are left just needing to show that $\F_p[y]$ has a nice filtration, when viewed as an object in $\F_p[y^{p^j}]-\F_p[y]-\U$.

By downwards induction on $i$, we show that, for $0 \leq i \leq j$, $\F_p[y^{p^i}]$ has a nice filtration, when viewed as an object in $\F_p[y^{p^j}]-\F_p[y]-\U$.  The case $i=j$ is clear.

For the inductive step, we filter $\F_p[y^{p^i}]$.  For $0 \leq r \leq p-1$, define $M(r)$ to be the span of $\{ y^{p^im} \ | \ m\equiv s \mod p, \text{ for some } 0 \leq s \leq r\}$.  Using the formulae $\PP^ky^n = \binom{n}{k}y^{n+k(p-1)}$ and $\Delta(y^n) = \sum_k \binom{n}{k} y^k \otimes y^{n-k}$, one easily checks that each $M(r)$ is an object in $\F_p[y^{p^j}]-\F_p[y]-\U$: if $\binom{n}{k} \not \equiv 0 \mod p$, and $n$ has the form $p^i(pa+s)$ with $0 \leq s \leq r \leq p-1$, then both $n+k(p-1)$ and $k$ also have this form.

Thus we have a filtration in $\F_p[y^{p^j}]-\F_p[y]-\U$:
$$ \F_p[y^{p^{i+1}}] = M(0) \subseteq M(1) \subseteq \dots \subseteq M(p-1) = \F_p[y^{p^{i}}],$$
and we are assuming by induction that $ \F_p[y^{p^{i+1}}]$ has a good filtration.  Now one checks that $M(r)/M(r-1) \simeq \Sigma^{2p^ir}M(0)$ as objects in $\F_p[y^{p^j}]-\F_p[y]-\U$, so by upwards induction on $r$ we conclude that each $M(r)$ has a good filtration.
\end{proof}

\begin{proof}[Proof of \lemref{d_1R_0, p odd}]  By \propref{d1 prop 1}(c), it suffices to prove that $d_1(\F_p[y]) = 1$, and, if $j>0$,
$ d_1(\F_p[y^{p^j}]) = 2p^{j-1}$.

\corref{d1 cor} (or \propref{d1 prop 2}(b)), applied to the short exact sequence
$$ 0 \ra \F_p[y] \ra \Lambda^*(x) \otimes \F_p[y] \ra \Sigma \F_p[y] \ra 0,$$
shows that $d_1(\F_p[y]) = d_0(\Sigma \F_p[y]) = 1$.

If $j > 1$, we consider the short exact sequence:
$$ 0 \ra \F_p[y^{p^j}] \ra \F_p[y^{p^{j-1}}] \ra  \F_p[y^{p^{j-1}}]/\F_p[y^{p^{j}}]\ra 0.$$
We claim that $d_0$ of the last term is $2p^{j-1}$, which, by induction, will be strictly more than $d_1$ of the middle term.  Thus \propref{d1 prop 2}(b) applies to say that $ d_1(\F_p[y^{p^j}]) = 2p^{j-1}.$

To verify the claim, one checks that the map $\F_p[y^{p^{j-1}}] \ra \Sigma^{2p^{j-1}}\F_p[y^{p^{j-1}}]$ sending $y^{p^{j-1}n}$ to the ${2p^{j-1}}$th suspension of $ny^{p^{j-1}(n-1)}$ is a map of unstable $\A$--modules, and thus induces an embedding $\F_p[y^{p^{j-1}}]/\F_p[y^{p^{j}}] \hra \Sigma^{2p^{j-1}}\F_p[y^{p^{j-1}}]$ in $\U$. Since the range of this embedding is the ${2p^{j-1}}$th suspension of a reduced module, the same is true of the domain, which thus has $d_0 = {2p^{j-1}}$.
\end{proof}

\subsection{ $d_1(G)$ when $G$ has a $p$--central $p$--Sylow subgroup}  Now suppose that $G$ is not necessarily $p$--central, but has a $p$--central $p$--Sylow subgroup $P$. Here we show that then $d_1(G) = d_1(P)$.  As $d_1(G) \leq d_1(P)$ is always true, the point is to show that $d_1(G)$ is as big as it could be.

Let $e_{\omega} \in \F_p[Out(P)]$ be an idempotent chosen so that $\F_p[Out(P)]e_{\omega}$ is the projective cover of the one dimensional module $\omega(P)^{-1}$.

\begin{lem} $d_1(G) = d_1(P)$ if and only if $d_1(R_{e(P)}H^*(G)) = h(P)$, and either of these equalities are implied by $d_1(e_{\omega}R_0H^*(P)) = h(P)$.
\end{lem}

\begin{proof} As we proved that $d_1(P)=e(P)+h(P)$, we showed that $d_1(nil_{e(P)}H^*(P)) = e(P)+h(P)$ and $d_1(H^*(P)/nil_{e(P)}H^*(P)) < e(P)+h(P)$.  This second fact implies that $d_1(H^*(G)/nil_{e(P)}H^*(G)) < e(P)+h(P)$ also holds, as $H^*(G)$ is a direct summand of $H^*(P)$ in $\U$.
We conclude that $d_1(G) = d_1(P)$ if and only if $d_1(nil_{e(P)}H^*(G)) = e(P)+h(P)$.  As $d_1(nil_{e(P)}H^*(G)) = e(P) + d_1(R_{e(P)}H^*(G))$, we deduce that $d_1(G) = d_1(P)$ if and only if $d_1(R_{e(P)}H^*(G)) = h(P)$.

Now reasoning as in \secref{d0(P) subsection}, this last equality would follow if one could show that $d_1(e_{\omega}R_0H^*(P)) = h(P)$.
\end{proof}

We now sketch a proof that $d_1(e_{\omega}R_0H^*(P)) = h(P)$.  This involves redoing the calculation that $d_1(R_0H^*(P)) = h(P)$ in a way that allows one to keep track of the $Out(P)$--action.  Let $C=C(P)$.

The 0--line of the spectral sequence associated to $C \ra P \ra P/C$ is natural with respect to the action of $Out(P)$.  Thus the filtration studied in \secref{E0 section},
$$ H^*(C) = E_2^{0,*} \supset E_3^{0,*} \supset E_{2p+1}^{0,*} \supset \dots \supset E_{2p^{k-1}+1}^{0,*} \supset E_{2p^k+1}^{0,*} = E_{\infty}^{0,*} = R_0H^*(P),$$
is a filtration by unstable modules with an $Out(P)$--action.

Our work above shows that $$d_1(E_r^{0,*}) =
\begin{cases}
2p^{j-1} & \text{if } r = 2p^j+1 \text{ with } j \geq 1 \\
1 & \text{if } r = 3 \\ 0 & \text{if } r=2.
\end{cases}$$

We now suppose that $p>2$ and $k \geq 1$: the cases when $p=2$ or when $E_{\infty}^{0,*}$ equals $E_2^{0,*}$ or $E_3^{0,*}$ are similar and easier.  Recall that then $h(P) = 2p^{k-1}$. Using \propref{d1 prop 2} in the usual way, we conclude that $d_1(e_{\omega}R_0H^*(P)) = h(P)$ if and only if $d_0(e_{\omega}B) = 2p^{k-1}$, where $B = E_{2p^{k-1}+1}^{0,*}/E_{2p^k+1}^{0,*}$.

From \secref{E0 section}, we see that
\begin{multline} \label{nasty formula}
B = \Lambda^*(C_0^{\#})\otimes  S^*(\beta(C_1^{\#}) + \Phi\beta(C_2^{\#}) + \dots + \Phi^{k-1}\beta(C_{k}^{\#})) \\
\otimes S^*(\Phi^{k-1}\beta(C^{\#}/C_{k}^{\#}))/ S^*(\Phi^{k}\beta(C^{\#}/C_{k}^{\#})),
\end{multline}
where $ C_0^{\#} \subseteq C_1^{\#} \subseteq \dots \subseteq C_k^{\#} \subseteq C^{\#}$ is a filtration of $C^{\#}$ as an $Out(P)$--module.

As an unstable module, $B$ thus has the form
$$ M \otimes (S^*(\Phi^{k-1}\beta(V))/S^*(\Phi^{k}\beta(V))),$$
where $M$ is reduced.  Now one observes that $S^*(\beta(V))/S^*(\Phi\beta(V)) = \Sigma^2 N$ where $N$ is reduced.  Thus
\begin{equation*}
\begin{split}
S^*(\Phi^{k-1}\beta(V))/S^*(\Phi^{k}\beta(V)) &
= \Phi^{k-1}(S^*(\beta(V))/S^*(\Phi\beta(V))) \\
  & = \Phi^{k-1}(\Sigma^2 N) \\
  & = \Sigma^{2p^{k-1}}(\Phi^{k-1}N),
\end{split}
\end{equation*}
which is the $2p^{k-1}$st suspension of a reduced module.

We conclude that $d_0(e_{\omega}B) = 2p^{k-1}$ if and only if $e_{\omega}B$ is nonzero.

The image of $Out(P) \ra GL(C)$ lands in the parabolic subgroup $GL(C,P)$ respecting the filtration of $C$, and the idempotent $e_{\omega}$ will project to a nonzero idempotent in $\F_p[GL(C,P)]$.

We claim that if $e \in \F_p[GL(C,P)]$ is any nonzero idempotent, then $e$ acts nontrivially on $B$ as described in (\ref{nasty formula}). Equivalently, we claim that all irreducible $\F_p[GL(C,P)]$--modules occur as composition factors in $B$.

To prove the claim, we note that all irreducible $GL(C,P)$ modules will be pullbacks from the associated Levi factor (i.e. the product of `block diagonal' $GL(V_j)$'s), as the projection from the one to the other has kernel which is a $p$--group.  This reduces us quickly to verifying the following lemma.

\begin{lem}  Every irreducible $\F_p[GL(V)]$--module occurs as a composition factor in $S^*(\beta(V))/S^*(\Phi \beta(V))$.
\end{lem}

\begin{proof}  It is well known that every such irreducible $S$ occurs in $S^*(\beta(V))$. Choosing an occurrence of lowest polynomial degree, it is clear that it will remain nonzero in the quotient $S^*(\beta(V))/S^*(\Phi \beta(V))$.
\end{proof}

\section{Central essential cohomology} \label{cess section}  Recall that $Cess^*(G)$ is defined to be the kernel of the restriction map
$$ H^*(G) \ra \prod_{\substack{C(G)< U \\ C(G) \neq U}} H^*(C_G(U)).$$
The invariants $e^{\prime}(G)$ and $e^{\prime \prime}(G)$ are then defined by letting
$$ e^{\prime}(G) = \max\ \{ d \ | \ Q_{A}Cess^d(G) \neq 0 \} \cup \{-1\}, $$
where $A$ is a Duflot subalgebra of $H^*(G)$, and
$$ e^{\prime \prime}(G) = \max\ \{ d \ | \ P_{C(G)}Cess^d(G) \neq 0 \} \cup \{-1\}. $$

In this section we study $Cess^*(G)$, $e^{\prime}(G)$, and $ e^{\prime \prime}(G)$, and connect them to the invariant $d_0(G)$.

\subsection{The structure of $Cess^*(G)$}  We begin by proving \thmref{Cess thm}. Most of this theorem is restated in the following.  We let $C=C(G)$, as usual.

\begin{prop} \label{Cess prop}  If $A$ is a Duflot subalgebra of $H^*(G)$, then the following hold. \\

\noindent{\bf (a)}  $Cess^*(G)$ is a free $A$--module. \\

\noindent{\bf (b)} The composite $P_C Cess^*(G) \hra Cess^*(G) \era Q_ACess^*(G)$ is monic. \\

\noindent{\bf (c)} The sequence $\displaystyle 0 \ra Q_ACess^*(G) \ra Q_AH^*(G) \ra \prod_{C(G)<U} Q_AH^*(C_G(U))$
is exact.
\end{prop}
\begin{proof}  It is convenient to let $M_1 = Cess^*(G)$, $M_2 = H^*(G)$, $M_3 = M_2/M_1$, $\displaystyle M_4 = \prod_{C(G)<U} H^*(C_G(U))$, and $M_5 = M_4/M_3$.  We have short exact sequences of $A$--modules
$$ 0 \ra M_1 \ra M_2 \ra M_3 \ra 0,$$
and 
$$ 0 \ra M_3 \ra M_4 \ra M_5 \ra 0.$$
Statements (a) and (c) of the proposition follow if we verify that all the $M_j$ are free $A$--modules. This we do by arguing as in \cite[proof of Thm.12.3.3]{carlson et al}.  

The homomorphisms $C \times C_G(V) \ra C_G(V)$ make each $M_j$ into a $H^*(C)$ comodule and each $M_j \otimes H^*(C)$ into an $A$--module such that the comodule structure map $M_j \ra M_j \otimes H^*(C)$ is a map of $A$--modules, and the composite $ M_j \ra M_j \otimes H^*(C) \ra M_j$ is the identity.  Thus $M_j$ is a direct summand of $M_j \otimes H^*(C)$ as an $A$--module, and we conclude that $M_j$ is free if $M_j \otimes H^*(C)$ is.

To show that $M_j \otimes H^*(C)$ is a free $A$--module, we give it a decreasing $A$--module filtration by letting $F^n = M_j^{\geq n} \otimes H^*(C)$.  Then each $F^n/F^{n+1} = M_j^n \otimes H^*(C)$ is a direct sum of copies of $H^*(C)$, and so is a free $A$--module, and the freeness of $M_j$ easily follows.

Finally, statement (b) follows from consideration of the commutative square
\begin{equation*}
\xymatrix{
 P_CCess^*(G) \ar[d] \ar[r] & Q_ACess^*(G) \ar[d]  \\
P_CH^*(G) \ar[r] & Q_AH^*(G). }
\end{equation*}
The left map is clearly monic, and \thmref{prim thm} says that the bottom map is also.  Thus so is the top map.
\end{proof}

To finish the proof of \thmref{Cess thm}, it remains to show that $Cess^*(G)$ is finitely generated as an $A$--module (with the corollary that $e^{\prime}(G)$ and $ e^{\prime \prime}(G)$ are well defined finite numbers).  As $A$ has Krull dimension equal to the rank of $C$, and $Cess^*(G)$ is a free $A$--module, it is equivalent to prove the next result.

\begin{prop}  The Krull dimension of $Cess^*(G)$ is at most the rank of $C$.
\end{prop}
\begin{proof}  The proposition follows from a result of Carlson.  It is convenient to let $R=H^*(G)$, and $I=Cess^*(G)$.  By definition, the Krull dimension of the $R$--module $I$ is the Krull dimension of the algebra $R/Ann(I)$.

Let $J$ be the image of
$$ \sum_{\substack{C < U \\ C \neq U}} \Ind_{C_G(U)}^G: \bigoplus_{\substack{C < U \\ C \neq U}} H^*(C_G(U)) \ra H^*(G).$$
By standard arguments\footnote{This follows immediately from the fact that $\Ind_{C_G(U)}^G: H^*(C_G(U)) \ra H^*(G)$ is a map of $H^*(G)$--modules.}, $J \subset Ann(I)$.  Thus $R/J \ra R/Ann(I)$ is a surjection, and so the Krull dimension of $R/Ann(I)$ is at most the Krull dimension of $R/J$.

In the notation of \cite{carlson depth}, $J = J_{c+1}$, where $c$ is the rank of $C$. Then \cite[Cor.2.2]{carlson depth} says that the Krull dimension of $R/J$ is at most $c$.
\end{proof}

The next proposition is easily verified.

\begin{prop} \label{cess prod prop} $Cess^*(G \times H)$ is naturally isomorphic to $Cess^*(G) \otimes Cess^*(H)$.  Thus $e^{\prime}(G\times H) = e^{\prime}(G)+ e^{\prime}(H)$ and $e^{\prime \prime}(G\times H) = e^{\prime \prime}(G)+ e^{\prime \prime}(H)$.
\end{prop}

We end this subsection by proving \propref{depth prop}, which we recall here.

\begin{prop} \label{depth prop 2} Assuming that $V = C(C_G(V))$, $Cess^*(C_G(V)) = \mathbf 0$ unless the rank of $V$ is at least equal to the depth of $H^*(G)$.
\end{prop}
\begin{proof}  Let $r(U)$ denote the $p$--rank of an elementary abelian $p$--group $U$, and let $d$ be the depth of $H^*(G)$.  Suppose that $V = C(C_G(V))$ and $r(V)<d$.  We wish to show that $Cess*(C_G(V)) = \mathbf 0$.  Note that, if $V < U$, then $C_{C_G(V)}(U) = C_G(U)$.  Thus $Cess^*(C_G(V))$ is the kernel of the restriction map
$$ H^*(C_G(V)) \ra \prod_{\substack{V < U < C_G(V) \\ V \neq U}} H^*(C_G(U)).$$
The kernel of this is contained in the kernel of $$ f: H^*(C_G(V)) \ra \prod_{\substack{V < U < C_G(V) \\ r(U)= d}} H^*(C_G(U)).$$
Thus it suffices to show that $f$ is monic.

Meanwhile, Carlson's theorem \cite{carlson depth} implies that the product of restriction maps
$$ g: H^*(G) \ra \prod_{\substack{U<G \\ r(U)=d}} H^*(C_G(U))$$
is monic.

There is a commutative diagram
\begin{equation*}
\xymatrix{
\coprod_{V < U < C_G(V), \  r(U)=d} V \times C_G(U) \ar[d] \ar[r] & V \times C_G(V) \ar[d]  \\
\coprod_{U < G, \ r(U)=d} C_G(U) \ar[r] & G. }
\end{equation*}
This induces a commutative diagram in $\U$
\begin{equation*}
\xymatrix{
H^*(V) \otimes H^*(C_G(V))  \ar[d] \ar[r]^-{1 \otimes f} &  \prod_{V < U < C_G(V), \  r(U)=d} H^*(V) \otimes H^*(C_G(U))\ar[d]  \\
H^*(G) \ar[r]^-{g} & \prod_{U < G, \ r(U)=d} H^*(C_G(U)). }
\end{equation*}
Adjointing, we get a commutative diagram
\begin{equation*}
\xymatrix{
H^*(C_G(V))  \ar[d] \ar[r]^-{f} &  \prod_{V < U < C_G(V), \  r(U)=d} H^*(C_G(U))\ar[d]  \\
T_VH^*(G) \ar[r]^-{T_Vg} & \prod_{U < G, \ r(U)=d} T_VH^*(C_G(U)). }
\end{equation*}

The left vertical map here is an inclusion, and the exactness of $T_V$ shows that $T_Vg$ is monic.  We conclude that $f$ is also monic.
\end{proof}

\subsection{Proof of \thmref{d0 thm}.}

\thmref{d0 thm} says that $$d_0(G) = \max\{ e^{\prime \prime}(C_G(V)) \ | \ V < G \}.$$
In the right hand side of this equation, one can restrict to subgroups $V$ such that $V = C(C_G(V))$, because if \ $ V < U$ and $U$ is central in $C_G(V)$, then $C_G(V)=C_G(U)$, so that $e^{\prime \prime}(C_G(V)) = e^{\prime \prime}(C_G(U))$.  Thus \thmref{d0 thm} will follow from the next theorem.

\begin{thm} \label{d0 thm2} $\bar R_dH^*(G) \neq 0$ if and only if  $P_VCess^d(C_G(V)) \neq 0$ for some $V < G$ satisfying $V = C(C_G(V))$.
\end{thm}

As we begin the proof of this, it is convenient to let, for $V < G$,
$$ Ess^*(V) = \ker \{ P_VH^*(C_G(V)) \ra \prod_{V < U} P_VH^*(C_G(U)) \},$$
with the product over all $U$ that are strictly bigger than $V$.  One easily verifies that
\begin{equation*}
Ess^*(V) =
\begin{cases}
P_VCess^d(C_G(V)) & \text{if } V = C(C_G(V)) \\ 0 & \text{otherwise.}
\end{cases}
\end{equation*}
Thus we wish to show that $\bar R_dH^*(G) = 0$ if and only if $Ess^d(V)=0$ for all $V <G$.

We prove the `only if' implication first.  Let $c_V \in \F_p[GL(V)]$ be the top Dickson invariant, and let $W_G(V) = N_G(V)/V$.  Then \cite[Lemma 7.8]{k4} (proved by using the formula in this paper's \propref{good Rd formula prop}) says that, for all $V$, there is an embedding
$$  (c_VH^*(V) \otimes Ess^d(V))^{W_G(V)} \subseteq \bar R_dH^*(G).$$
Thus $\bar R_dH^*(G) = 0$ implies that  $(c_VH^*(V) \otimes Ess^d(V))^{W_G(V)} = 0$ for all $V$.  But then $Ess^d(V) = 0$, by the next lemma.

\begin{lem} If $W < GL(V)$ and $M$ is an $\F_p[W]$--module, then $[c_VH^*(V) \otimes M]^W = 0$ implies that $M=0$.
\end{lem}

\begin{proof} It is well known (see, e.g. \cite[p.45]{alperin}) that $\F_p[W]$ embeds in $H^*(V)$ as $\F_p[W]$--modules, and thus in $c_VH^*(V)$, as multiplication by $c_V$ is a monic $GL(V)$--module map. Thus $(\F_p[W] \otimes M)^W$ embeds in $(c_VH^*(V) \otimes M)^W$.  But $\F_p[W] \otimes M \simeq \F_p[W] \otimes M_{triv}$ as $\F_p[W]$--modules, where $M_{triv}$ denotes $M$ with trivial $W$--action.  Finally, $(\F_p[W] \otimes M_{triv})^W \simeq M$ as $\F_p$--vector spaces.  Putting this all together, we have shown that $M$ embeds in $(c_VH^*(V) \otimes M)^W$, and the lemma follows.
\end{proof}

Returning to the proof of the theorem, we now assume that $Ess^d(V) = 0$ for all $V$, and deduce that $\bar R_dH^*(G) = 0$.  Using the formula in \propref{good Rd formula prop}, given
$$\displaystyle x = (x_V) \in \bar R_dH^*(G) \subseteq \prod_V H^*(V) \otimes P_VH^d(C_G(V)),$$ we show that each component $x_V$ of $x$ is zero by downwards induction on $V$.

So assume that $x_{U} = 0$ for all $V < U$.  \propref{good Rd formula prop} tells us that, under the restriction map
$$H^*(V) \otimes P_VH^d(C_G(V)) \ra \prod_{V < U} H^*(V) \otimes P_VH^d(C_G(U)),$$
$x_V$ will have the same image as $(x_U)$ under the map
$$ \prod_{V < U} H^*(U) \otimes P_UH^d(C_G(U)) \ra \prod_{V < U} H^*(V) \otimes P_VH^d(C_G(U)).$$
Since the latter is zero by inductive assumption, we conclude that $x_V \in H^*(V) \otimes Ess^d(V)$, and is thus zero also.

\begin{rem} Note that \thmref{d0 thm2} includes a second proof that $e^{\prime \prime}(G)$ is a well defined finite number.
\end{rem}

\subsection{The Depth Conjecture, the Regularity Conjecture, and a bound on $e^{\prime}(G)$.}

By \thmref{Cess thm}, $e^{\prime \prime}(G) \leq e^{\prime}(G)$.  Thus we get the bound
$$d_0(G) \leq \max\{ e^{\prime}(C_G(V)) \ | \ V < G \}.$$

If $G$ is $p$--central, we know that $e^{\prime \prime}(G) = e^{\prime}(G) = e(G)$.  Here we discuss work towards \conjref{e' conjecture} which said that, if $G$ is not $p$--central, then $e^{\prime}(G) < e(G)$.  In particular, we show that this is true if the $p$--rank of $G$ is no more than 2 more than the rank of $C(G)$, and link the general conjecture to Benson's Regularity Conjecture.  Enroute, a similar argument will also prove \thmref{Cess depth thm}, the special case of Carlson's Depth Conjecture in which the depth of $H^*(G)$ is as small as possible.

Thanks to \propref{cess prod prop}, it suffices to prove either \thmref{Cess depth thm} or \conjref{e' conjecture} in the special case when $G$ has no direct summands isomorphic to $\Z/p$.  We will assume this.  As a consequence, even in the odd prime case, a Duflot subalgebra $A$ of $H^*(G)$ will have the form
$$ A = \F_p[\xi_1, \dots, \xi_c],$$
where $c$ is the rank of $C= C(G)$ and $\displaystyle e(G) = \sum_{j=1}^c (|\xi_j| -1)$.
The sequence $\xi_1, \dots, \xi_c$ is a {\em Duflot sequence}: the sequence restricts to a regular sequence in $H^*(C)$. We let $r$ be the $p$--rank of $G$.

\begin{lem} Suppose $c<r$.  Given any Duflot sequence $\xi_1, \dots, \xi_c \in H^*(G)$, there exists $\xi \in H^*(G)$ such that, for all proper inclusions $C<V$, $\xi_1, \dots, \xi_c, \xi$ restricts to a regular sequence in $H^*(C_G(V))$.
\end{lem}

\begin{proof}  Let $n$ be the rank of $G/C$ (so $n > r-c$), and let $\rho$ be the regular representation of $G/C$.  For $1 \leq i \leq r-c$, let
$\bar \kappa_i \in H^{2(p^n-p^{n-i})}(G/C)$ be the $(p^n-p^{n-i})$th Chern class of $\rho$, and then let $\kappa_i = \Inf_{G/C}^G(\bar \kappa_i) \in H^*(G)$.  It is easy to check that $\xi_1, \dots, \xi_c, \kappa_1, \dots, \kappa_{r-c}$ is a polarized system of parameters in the sense of \cite[Def.2.2]{green1}.  It follows that the element $\xi = \kappa_1$ satisfies the conclusion of the lemma.
\end{proof}

\begin{prop}  \label{handy prop} Suppose $c<r$ and $A = \F_p[\xi_1, \dots, \xi_c]$ is a Duflot algebra of $H^*(G)$. The following are equivalent for a fixed integer $e \geq 0$. \\

\noindent{\bf (a)} $e^{\prime}(G) < e$. \\

\noindent{\bf (b)} With $\xi$ as in the lemma, the kernel of multiplication by $\xi$,
$$ \ker \{\xi \cdot : H^d(G)/(\xi_1, \dots \xi_c) \ra H^{d+|\xi|}(G)/(\xi_1, \dots \xi_c)\},$$
is zero for all $d \geq e$. \\

\noindent{\bf (c)} $\displaystyle \bigcap_{\xi \in \tilde H^*(G)} \ker \{\xi \cdot : H^d(G)/(\xi_1, \dots \xi_c) \ra H^{d+|\xi|}(G)/(\xi_1, \dots \xi_c)\}$ \\
is zero for all $d \geq e$.
\end{prop}

\begin{proof} For each $d$ and $\xi \in H^*(G)$, we have a commutative  diagram
\begin{equation} \label{key diagram}
\xymatrix{
H^d(G)/(\xi_1, \dots \xi_c) \ar[d]^{\xi \cdot} \ar[r]^-{f(d)} & \bigoplus_{C<V}H^d(C_G(V))/(\xi_1, \dots \xi_c) \ar[d]^{\xi \cdot}  \\
H^{d+|\xi|}(G)/(\xi_1, \dots \xi_c) \ar[r]^-{f(d+|\xi|)} & \bigoplus_{C<V}H^{d+|\xi|}(C_G(V))/(\xi_1, \dots \xi_c), }
\end{equation}
where $f(d)$ is induced by the evident restriction maps.

By \thmref{Cess thm}, $f(d)$ is monic for all large $d$, and $e^{\prime}(G)$ is the largest $d$ such that $f(d)$ is not monic.  Thus statement (a) is equivalent to the statement that $\ker f(d)$ is zero for all $d \geq e$.

We show that statement (a) implies statement (b).  Thus suppose that $\ker f(d)$ is zero for all $d \geq e$, and let $\xi \in H^*(G)$ be as in the lemma.  Then, in diagram (\ref{key diagram}), the right map is monic for all $d$, and the top map is monic for all $d \geq e$.  Thus the left map is monic in the same range.

Statement (b) obviously implies statement (c).

Finally, we show that statement (c) implies statement (a).  Assuming statement (c), we prove, by downwards induction on $d$, that $f(d)$ is monic.  Thus assume $f(d^{\prime})$ is monic for all $d^{\prime}>d \geq e$.  Given $0 \neq \kappa \in H^d(G)/(\xi_1, \dots \xi_c)$, we need to show that $f(d)(\kappa) \neq 0$.  By (c), there exists $\xi \in \tilde H^*(G)$ such that $\xi \cdot \kappa \neq 0$.  As $f(d+|\xi|)$ is monic by inductive assumption, $f(d+|\xi|)(\xi \cdot \kappa) \neq 0$.  But this equals $\xi \cdot f(d)(\kappa)$, and so $f(d)(\kappa) \neq 0$.
\end{proof}

\begin{proof}[Proof of \thmref{Cess depth thm}] With notation as in the proposition just proved, we wish to prove that $Cess^*(G) = 0$ if and only if the depth of $H^*(G)$ is greater than $c$. Thanks to \thmref{Cess thm}, $Cess^*(G)=0$ if and only if statement (a) of the last proposition holds when $e=0$. But then statement (b) is true with $e=0$, and thus the depth of $H^*(G)$ is at least $c+1$.

Conversely, if the depth of $H^*(G)$ is at least $c+1$, there exists a $\xi$ such that $\xi_1, \dots, \xi_c, \xi$ is a regular sequence on $H^*(G)$, and so statement (c) certainly holds with $e=0$. Thus statement (a) does as well.
\end{proof}

Now we study statements (b) and (c) of the last proposition, using work by Carlson and Benson.

\begin{prop} \label{r-c=1 prop} If \ $r-c = 1$, then $Q_ACess^*(G)$ satisfies Poincar\'e duality with duality degree equal to $e(G)$.  In other words, the Poincar\'e polynomial $p_{Q_ACess^*(G)}(t)$ satisfies
$$ p_{Q_ACess^*(G)}(t) = t^{e(G)}p_{Q_ACess^*(G)}(1/t).$$
\end{prop}

\begin{proof}  The conclusion of the proposition is obvious if $Cess^*(G) = 0$, so we can assume that the depth of $H^*(G)$ is precisely $c$.  Let $\xi_1, \dots, \xi_c$ be as in \propref{handy prop}, and choose $\xi$ as in the lemma.  Replacing $\xi$ by a large power of itself, if necessary, we can assume that, in diagram (\ref{key diagram}), $f(d+|\xi|)$ is monic for all $d$.  Thus $Q_ACess^*(G)$, the kernel of the top map in (\ref{key diagram}), identifies with the kernel of multiplication by $\xi$, the left map in (\ref{key diagram}).  But a careful reading of \cite[Lemma 3.2 and its proof]{bc2} reveals  that the Poincar\'e series of this kernel is precisely the polynomial called `$p_{r}(t)$' there, and then \cite[Theorem 3.9]{bc2} says that the functional equation of the proposition holds.
\end{proof}

\begin{cor} If \ $r-c = 1$, then
$$ e^{\prime}(G) = e(G) - \min \{ d \ | \ Cess^d(G) \neq 0\} < e(G).$$
\end{cor}

To state what we know about the situation when $r-c>1$, we need to introduce local cohomology.  If $I$ is a homogeneous ideal in a graded ring $R$, and $M$ is a graded $R$--module, $H^{0,*}_I(M)$ is defined to be the $I$--torsion in $M$, i.e. the set of $x \in M$ such that $I^kx = 0$ for some $k$.  This is a left exact functor of $M$, and $H^{d,*}_I(M)$ is defined to be the associated $d$th right derived functor.

In \cite{b}, Benson conjectured

\begin{conj}[Strong Regularity Conjecture]
\begin{equation*}
H_{\tilde H^*(G)}^{i,j}(H^*(G)) =0 \text{ for }
\begin{cases}
j \geq -i & \text{if } c \leq i < r \\ j > -i & \text{if } i=r.
\end{cases}
\end{equation*}
\end{conj}

\begin{proof}[Proof of \propref{conjectures prop}] This proposition asserted that, for a fixed finite group $G$, \conjref{e' conjecture} is implied by the Strong Regularity Conjecture.

In the terminology of \cite{b}, if the Strong Regularity Conjecture holds, then, by \cite[Thm.4.5]{b}, every filter regular sequence is of type beginning with the sequence $(-1, -2, \dots, -(c+1))$.  In particular, with $\xi$ as in statement (b) of \propref{handy prop}, the sequence $\xi_1, \dots, \xi_c, \xi$ is the beginning of such a sequence. From the definition of filter regular, we see that statement (b) of \propref{handy prop} thus holds with $e=e(G)$.
\end{proof}

As mentioned in the introduction, in \cite{b}, Benson shows that his conjecture is true if $r-c \leq 2$.  I have my own `heuristic' proof of statement (c) with $e=e(G)$ under the same condition, and the failure of the method to go beyond $r-c \leq 2$ makes one wonder if a counterexample to both of our conjectures is lurking among the groups of order 128 or 256.

\begin{rem} Slightly milder than the Strong Regularity Conjecture is Benson's Regularity Conjecture, which asserts that the Castelnuovo--Mumford regularity of $H^*(G)$ is precisely 0.  In terms of local cohomology, this is the statement that
$H_{\tilde H^*(G)}^{i,j}(H^*(G)) =0$ for $j>-i$.  For our purposes, this is enough to deduce that $e^{\prime}(G) \leq e(G)$.
\end{rem}

\section{Examples} \label{ex section}

\begin{ex} \label{W(2) ex}  Let $W(2)$ be the universal 2--central group whose quotient by its center $C$ is $V_2=(Z/2)^2$.  Thus there is a central extension
$$ H_2(V_2;\F_2) \xra{i} W(2) \xra{q} V_2,$$
where $C=H_2(V_2;\F_2) \simeq (\Z/2)^3$.  In terms of Hall--Senior numbering, and thus also the numbering in \cite{carlson et al}, $W(2)$ is 32\#18.

In the associated spectral sequence, one has that
$$ E_2^{*,*} = \F_2[x,y,a,b,c],$$
with $a,b,c \in E_2^{0,1}$ and $x,y \in E_2^{1,0}$, and $d_2(a)=x^2$, $d_2(b)=xy$, and $d_2(c)=y^2$.

As $E_3^{*,0} = \F_2[x,y]/(x^2,xy,y^2)$, it follows that $a^2,b^2$, and $c^2$ must be permanent cycles.  We conclude that $W(2)$ will have type $[2,2,2]$ so that $d_0(W(2)) = e(W(2)) = 3$ and $d_1(W(2)) = 4$.

With a bit more work, one can show that $E_3^{*,*}/(a^2,b^2,c^2)$ is six dimensional with generators as indicated in Figure \ref{figure 1}, where $u$ and $v$ are respectively represented by $bx+ay$ and $cx+by$.  This is a Poincar\'e duality algebra with relations $x^2 = xy = y^2 = u^2=v^2=uv=xu = yv = xv+yu=0$.
\begin{figure}
\begin{equation*}
\begin{array}{c|ccccc}
q &&&& \\
  \\
2 &&&&& \\
1 &&u,v& xv && \\
0 &1&x,y&&& \\ \hline
 &0&1&2& &p \\
\end{array}
\end{equation*}
\caption{$E_3^{p,q} = E_{\infty}^{p,q}$ modulo $(a^2,b^2,c^2)$} \label{figure 1}
\end{figure}

It follows that $E_3^{*,*} = E_{\infty}^{*,*}$, and then that
$$ H^*(W(2)) \simeq \F_2[\alpha, \beta, \gamma, x,y,u,v]/(x^2,xy,y^2,u^2,uv,v^2,xu,yv, xv+yu),$$
with $x,y \in H^1$ and  $\alpha,\beta,\gamma,u,v \in H^2$. Here $\alpha$, $\beta$, and $\gamma$ are represented by $a^2$, $b^2$, and $c^2$ in the spectral sequence.

The polynomial subalgebra $A= \F_2[\alpha, \beta, \gamma]$ is a Duflot subalgebra.  With respect to the $H^*(C) = \F_2[a,b,c]$  comodule structure, the elements $1,x,y$ are in the image of the inflation map $q^*$ and so are primitive. The top class $xv$ is not in the image of inflation, but is primitive, by our general theory. The elements $u$ and $v$ are {\em not} primitive, as $m^*(u) = 1 \otimes u + b \otimes x + a \otimes y$ and $m^*(v) = 1 \otimes v + c \otimes x + b \otimes y$ in $H^*(C) \otimes E_{\infty}^{*,*}$.  Thus each of the inclusions
$$ \im(q^*) \hra P_CH^*(W(2)) \hra Q_AH^*(W(2))$$
is {\em proper}.

The nilpotent filtration works as follows.

$R_0 = H^*(W(2))/(x,y,u,v) \simeq \F_2[\alpha, \beta, \gamma]$ and $\bar R_0 = \F_2[a,b,c]$.  The embedding $R_0 \subset \bar R_0$ sends $\alpha$ to $a^2$, $\beta$ to $b^2$, and $\gamma$ to $c^2$.

$R_1$ is the free $\F_2[\alpha, \beta, \gamma]$--module on generators $\bar x, \bar y, \bar u, \bar v$ of respective degrees 0,0,1,1, and $\bar R_1$ is the free $\F_2[a,b,c]$--module on $\bar x, \bar y$.  The embedding $R_1 \subset \bar R_1$ sends $\bar u$ to $b\bar x + a \bar y$ and $\bar v$ to $c\bar x + b\bar y$. Thus  $Sq^1(\bar u) = \beta \bar x + \alpha \bar y$ and $Sq^1(\bar v) = \gamma \bar x + \beta \bar y$.

$R_2$ and $\bar R_2$ are both 0, as $P_CH^2(W(2)) = 0$.

$R_3$ is the free $\F_2[\alpha, \beta, \gamma]$--module on a single generator $\bar x \bar v$ of degree 0, and $\bar R_3$ is the free $\F_2[a,b,c]$--module on this same element.

Finally, $H^*(W(2))_{LF} \subset H^*(W(2))$ is the algebra spanned by $1,x,y,xv$.  All nontrivial products and Steenrod operations are zero.
\end{ex}

\begin{ex} \label{group 64 no.108} Let $G$ be the group of order 64 with Hall--Senior number \#108. Using information from \cite{carlson et al}, we analyzed $H^*(G)$ in detail for other purposes in \cite{k4}.  Here we summarize relevant bits to illustrate how one can calculate $H^*(G)_{LF}$ and $\bar R_dH^*(G)$ by using \propref{good LF formula prop} and \propref{good Rd formula prop}.

The commutator subgroup $Z=[G,G]$ has order 2. The center $C$ is elementary abelian of rank 2, and $C = \Phi(G)$, so $Z < C$ and $G/C$ is elementary abelian of rank 4.   There is a unique maximal elementary abelian group $V$ of rank 3, and its centralizer $K$ has order 32, so that $N_G(V)/C_G(V) = G/K \simeq \Z/2$.  More precisely, $K$ is isomorphic to $(\Z/2)^2  \times Q_8$, with $Q_8$ embedded so that $V \cap Q_8 = Z$.

We have the following picture of $\A^C(G)$:
\begin{equation*}
\xymatrix{
C  \ar[r] & V \ar@(ur,dr)^{\Z/2}   }
\end{equation*}
and from this it is already clear that $\bar R_0H^*(G) = H^*(V)^{\Z/2}$.

We have maps of unstable algebras equipped with $Aut(G)$ action:
\begin{equation*} P_VH^*(K)\hookrightarrow P_CH^*(K) \xla{j^*} P_CH^*(G),
\end{equation*}
where $j:K \ra G$ is the inclusion. It is easily checked that
$j^*$ is onto in degree 1.

The maps of pairs $(Q_8,Z) \ra ((\Z/2)^2 \times Q_8, (\Z/2)^2 \times Z) = (K,V)$ induces an isomorphism of algebras:
$$P_VH^*(K) \simeq P_ZH^*(Q_8).$$

The algebra $P_ZH^*(Q_8)$ is familiar: the calculation of $H^*(Q_8)$ using the Serre spectral sequence associated to $Z \ra Q_8 \ra Q_8/Z$ reveals that $P_ZH^*(Q_8) = \text{Im }\{ H^*(Q_8/Z) \ra H^*(Q_8)\} = B^*$, where $B^*$ is the Poincar\'e duality algebra $\F_2[x,w]/(x^2 + xw + w^2, x^2w + xw^2)$, where $x$ and $w$ both have degree 1.  $B^*$ has dimension 1,2,2,1 in degrees 0,1,2,3.

From this we learn that $P_CH^*(K) \simeq B^*[y]$ where $y$ is also in degree 1, and thus is generated by elements in degree 1.  It follows that $j^*: P_CH^*(G) \ra P_CH^*(K)$ is onto, and then that $Inn(G)$ acts trivially on both $P_VH^*(K)$ and $P_CH^*(K)$.

\propref{good LF formula prop} then tells us that there is a pullback diagram of unstable algebras:
\begin{equation*}
\xymatrix{
H^*(G)_{LF} \ar[d] \ar[r] & P_CH^*(G) \ar[d]^{j^*}  \\
P_VH^d(K) \ar[r] & P_CH^*(K). }
\end{equation*}

Similarly, \propref{good Rd formula prop} tells us that, for all $d$, there is a pullback diagram of unstable modules:
\begin{equation*}
\xymatrix{
\bar R_dH^*(G) \ar[d] \ar[r] & H^*(C) \otimes P_CH^d(G) \ar[d]^{1\otimes j^*}  \\
H^*(V)^{\Z/2} \otimes P_VH^d(K) \ar[r] & H^*(C) \otimes P_CH^d(K). }
\end{equation*}

Note that the kernel of $j^*:H^*(G) \ra H^*(K)$ is precisely $Cess^*(G)$, which is described in \cite{carlson et al}.  In our terminology, we learn that a Duflot subalgebra $A$ is polynomial on classes of degree 2 and 8 (so $G$ has type $[8,2]$), and $Q_ACess^*(G)$ is a graded vector space of dimension 1,3,5,6,5,3,1 in degrees 1,2,3,4,5,6,7.  Note that this evident Poincar\'e duality is predicted by \propref{r-c=1 prop}.

$Q_ACess^*(G)$ has a basis in which every element is a product of 1 dimensional classes, and thus $P_CCess^*(G) = Q_ACess^*(G)$.  In \cite[Prop.10.2]{k4}, we further showed that
$$ P_CCess^*(G) \simeq \Sigma B^*[y]/(y^4)$$
as unstable modules.

It follows that there are short exact sequences in $\U$:
$$ 0 \ra \Sigma B^*[y]/(y^4) \ra H^*(G)_{LF} \ra B^* \ra 0,$$
and
$$ 0 \ra H^*(C) \otimes [\Sigma B^*[y]/(y^4)]^d \ra \bar R_dH^*(G) \ra H^*(V)^{\Z/2} \otimes B^d \ra 0.$$
Furthermore, $d_0(G) = e^{\prime \prime}(G) = e^{\prime}(G) =7$, and $e(G) = 8$.
\end{ex}

\begin{ex}  One can often determine $e(G)$ using minimal information about the extension class $\tau^*: C^* \ra H^2(G/C)$ (where $C = C(G)$), and in situations where $H^*(G)$ has yet to be calculated.

For example, suppose that $p=2$ and $G$ has no $\Z/2$ direct summands (so that $\tau^*$ is monic).  If the image of $\tau^*$ has a basis consisting of products of 1 dimensional classes, then $G$ has type $[2,\dots,2]$ and so $e(G)$ equals the rank of $C$.  To see this, we note that, if $d_2(a) = xy$, then
$$d_3(a^2) = d_3(Sq^1a) = Sq^1(d_2(a)) = Sq^1(xy) = x^2y+xy^2 \equiv 0 \mod (xy).$$

This criterion holds for the important family of groups studied in \cite{akm}.  There the authors associate a 2-central Galois group  $\mathcal G_{\F}$ to every field $\F$ of characteristic different from $2$ that is not formally real. They call this group a $W$--group due to its connections to the Witt ring $W\F$ \cite{mis}.  Thus $d_0(\mathcal G_{\F}) = e(\mathcal G_{\F}) = r$ and $d_1(\mathcal G_{\F}) = r+1$, where $\mathcal G_{\F}$ has rank $r$.  Included among these groups are the universal $W$--groups $W(n)$, the 2--central group with extension sequence
$$ H_2((\Z/2)^n;\F_2) \ra W(n) \ra (\Z/2)^n.$$
Thus $d_0(W(n)) = \binom{n+1}{2}$ and $d_1(W(n)) = \binom{n+1}{2}+1$.

At odd primes $p$, analogous criteria exist, ensuring that $G$ is $p$--central of type $[2,\dots,2]$.  Interesting families of such groups were studied by Browder--Pakianathan \cite{bp} and Adem--Pakianathan \cite{ap}.  Included among these are the universal groups $W(n,p)$,  with extension sequence
$$ H_2((\Z/p)^n;\F_p) \ra W(n,p) \ra (\Z/p)^n.$$
Thus $d_0(W(n,p)) = \binom{n+1}{2}$ and $d_1(W(n,p)) = \binom{n+1}{2}+1$.

\end{ex}

\begin{ex} Compared to the families in the last example, at the other extreme among $2$--central 2--groups is the 2--Sylow subgroup $P$ of the simple group $SU(3,4)$.  This group has order 64 and Hall--Senior number \#187.  Its center $C$ is elementary abelian of rank 2.  In \cite{green2}, Green analyzed the associated spectral sequence\footnote{A key simplification comes by computing with $\F_4$ coefficients rather than $\F_2$ coefficients.}.  In particular, $P$ has type $[8,8]$, so that $d_0(P) = e(P)=14$ and $d_1(P) = 18$, and the analogue of Figure \ref{figure 1} is the impressively complex Figure \ref{figure 2} (reproduced from \cite{green2}).
\begin{figure}
\begin{equation*}
\begin{array}{c|ccccccccccc}
q &&&& \\
  \\
8 &&&&2 \\
7 \\
6 & & & &6&8&8&8&4&1 \\
5 \\
4 & & &8&12&8&7&4 \\
3 \\
2 &&&4&7&8&8&8 \\
1 & \\
0 &1&4&8&10&8&6 \\ \hline
 &0&1&2&3&4&5&6&7&8& &p \\
\end{array}
\end{equation*}
\caption{The dimension of $E_{\infty}^{p,q}$ modulo $(a^8,b^8)$} \label{figure 2}
\end{figure}

In spite of this complexity, it is interesting to note that one can get the bound $e(P) \leq 14$ quite easily, by using representation theory and characteristic classes.

We thank David Green for the following description of some complex representations of $P$.  Let $H^*(C) = \F_2[a,b]$, and then let $\rho_a$ and $\rho_b$ be the 1--dimensional complex representations of $C$ with respective total Stiefel--Whitney classes $w(\rho_a) = 1 + a^2$, $w(\rho_b) = 1+ b^2$. These representations extend to 1-dimensional representations $\tilde \rho_a$ and $\tilde \rho_b$ of subgroups $Q_a$ and $Q_b$ of index 4 in $P$.  Let $\omega_a$ and $\omega_b$ be the 4 dimensional representations one gets by inducing
$\tilde \rho_a$ and $\tilde \rho_b$ up to $P$: these turn out to be irreducible.

By construction $\Res_C^P(\omega_a) = 4\rho_a$ and $\Res_C^P(\omega_b) = 4\rho_b$.  It follows that the total Stiefel--Whitney classes of $\omega_a$ and $\omega_b$ restrict to $(1+a^2)^4 = 1 + a^8$ and  $(1+b^2)^4 = 1+ b^8$ in $H^*(C)$.  Thus $\im \Res_C^P$ contains $\F_2[a^8,b^8]$ and so $e(P) \leq 14$ must hold.

Alternatively, one can just use the single 8 dimensional representation $\omega_a \oplus \omega_b$.  This is faithful, as it is faithful when restricted to $C$, the subgroup of all elements of order 2.  It has characteristic classes that restrict to $a^8 + b^8$ and $a^8b^8$ in $H^*(C)$.  From this, one can formally deduce that the special Hopf algebra $\im \Res_C^P$ must contain $\F_2[a^8,b^8]$, so that $d_0(P) \leq 14$ and $d_1(P) \leq 18$.  By contrast, the estimate of Henn, Lannes, and Schwartz in \cite{hls1} just lets one conclude that  $d_0(P) \leq 64$ and $d_1(P) \leq 120$ if one knows that $P$ has a faithful 8 dimensional complex representation.  This suggests that there might be some general bounds for $d_s(G)$ for an arbitrary group $G$, determined by the dimensions of its faithful representations, that are much better than those in \cite{hls1}.
\end{ex}

\appendix
\newpage

\section{Tables of group invariants}

Here are various tables of some of our invariants for 2--groups of order dividing 64. The tables were compiled by hand using the calculations in \cite{carlson et al} and the website version \cite{carlson website}.  The type of a group $G$, and thus $e(G)$ and $h(G)$, can be deduced by inspecting the description of restriction to maximal elementary abelian subgroups; this is particularly easy when $G$ is 2--central.  If $G$ is not 2--central, one can immediately determine if $Cess^*(G) \neq 0$, since both the rank of $Z(G)$ and the depth of $H^*(G)$ are given, and then read off the number $e^{\prime}(G)$ from the description of depth essential cohomology.  The website source allows one to identify centralizers of elementary abelian subgroups as needed.

We say a group is indecomposable if it cannot be written as a nontrivial direct product of two subgroups.  The numbering of groups is as in \cite{carlson et al} which follows the Hall--Senior numbering \cite{hall senior}.

In Tables 1 and 2, recall that, since $G$ is $2$--central, $d_0(G) = e(G) = e^{\prime}(G) = e^{\prime \prime}(G)$, and $d_1(G) = e(G)+h(G)$.

In Table 3, `$\mathbf 2$' means $\Z/2$, etc.  To compute $d_0(G)$, we needed to observe that, in all cases covered by this table, $e^{\prime \prime}(G) = e^{\prime}(G)$.  Except when $G$ is $32\#41$, this can be checked by noticing that elements in the top degree in $Q_ACess^*(G)$ are represented by classes in the image of $\Inf_{G/C}^G$, and so are primitive.  When $G$ is $32\#41$, elements in the $Q_ACess^5(G)$ are represented by essential classes of lowest degree, and so are primitive. \\

\begin{center}
{\bf Table 1: Indecomposable, 2--central, 2--groups of order $\leq 32$ }
\bigskip

\begin{tabular}{|c|c|c|c|c|c|} \hline
 & & & & & \\[-.12in]
Order & \# & Type  & $d_0(G)$ & $d_1(G)$ & Notes \\   \hline
 & & & & & \\[-.12in]
2 & 1 & [1] & 0 & 0 & $\Z/2$ \\   \hline
 & & & & & \\[-.12in]
4  & 2  & [2]  & 1  & 2  &  $\Z/4$     \\   \hline
   & & & & & \\[-.12in]
8  & 3  & [2]   & 1  & 2  &  $\Z/8$     \\   \hline
   & & & & & \\[-.12in]
  & 5  & [4]    & 3  & 5  &  $Q_8$     \\   \hline
   & & & & & \\[-.12in]
16  & 5  &  [2]   & 1  & 2  & $\Z/16$      \\   \hline
   & & & & & \\[-.12in]
  & 14  & [4]    & 3  & 5  & $Q_{16}$       \\   \hline
   & & & & & \\[-.12in]
32  & 18  & [2,2,2]    & 3   & 4  &       \\   \hline
   & & & & & \\[-.12in]
  & 19  &  [2,2]   &  2 & 3  &       \\   \hline
   & & & & & \\[-.12in]
  & 21  & [2,2]    &  2 &  3 &       \\   \hline
   & & & & & \\[-.12in]
  & 28  & [4,2]    & 4  &  6 &       \\   \hline
   & & & & & \\[-.12in]
  & 29  &  [2,2]   & 2  & 3  &       \\   \hline
   & & & & & \\[-.12in]
  & 30  & [2,2]    & 2  & 3  &       \\   \hline
   & & & & & \\[-.12in]
  & 35  & [4,2]    & 4 & 6  &       \\   \hline
   & & & & & \\[-.12in]
  & 40  &  [4,4]   & 6  & 8  &       \\   \hline
   & & & & & \\[-.12in]
  & 51  &  [4]   & 3  & 5  &  $Q_{32}$     \\   \hline
\end{tabular}
\end{center}
\newpage
\begin{center}
{\bf Table 2: Indecomposable, 2--central, groups of order 64 }
\bigskip

\begin{tabular}{|c|c|c|c|c|} \hline
  & & & & \\[-.12in]
 \# & Type  & $d_0(G)$ & $d_1(G)$ & Notes \\   \hline

    & & & & \\[-.12in]
 11  & [2]    & 1  & 2  &  $\Z/64$     \\   \hline

   30  & [2,2,2]    & 3  &  4 &       \\   \hline

   37  &  [2,2,2]   & 3  & 4  &       \\   \hline

   38  & [2,2]    & 2  & 3  &       \\   \hline

   39  & [2,2]    & 2  & 3  &       \\   \hline

   41  & [2,2]    & 2  & 3  &       \\   \hline

   59  & [2,2,2]    & 3  & 4  &       \\   \hline

    63 & [4,2]    & 4  & 6   &       \\   \hline

   64  & [2,2]    & 2  & 3  &       \\   \hline

   65  & [2,2]    & 2  & 3  &       \\   \hline

   82  & [2,2,2]    & 3  & 4  &       \\   \hline

   87  & [4,2,2]    & 5  & 7  &       \\   \hline

   88  & [2,2,2]    & 3  & 4  &       \\   \hline

   90  & [2,2,2]    & 3  & 4  &       \\   \hline

   92  & [4,2,2]    & 5  & 7  &       \\   \hline

   93  & [2,2,2]    & 3  & 4  &       \\   \hline

   101  & [4,4]    & 6  & 8  &       \\   \hline

   119  & [4,2]    & 4  & 6  &       \\   \hline

   139  & [4,2]    & 4  & 6  &       \\   \hline

   140  & [2,2]    & 2  & 3  &       \\   \hline

   141  & [2,2]    & 2  & 3  &       \\   \hline

   145  & [4,2,2]    & 5  & 7  &       \\   \hline

   149  & [4,2,2]    & 5  & 7  &       \\   \hline

   152  & [4,2,2]    & 5  & 7  &       \\   \hline
     & & & & \\[-.12in]
   153  & [4,4,4]    & 9  & 11  & 2--Sylow of $Sz(8)$      \\   \hline

   162  & [4,4]    & 6  & 8  &       \\   \hline
     & & & & \\[-.12in]
   187  & [8,8]    & 14  & 18  & 2--Sylow of $U_3(\F_4)$      \\   \hline

   190  & [4,2]    & 4  & 6  &       \\   \hline

   191  & [4,4]    & 6  & 8  &       \\   \hline

   192  & [4,2]    & 4  & 6  &       \\   \hline

   194  & [4,4]    & 6  & 8  &       \\   \hline

   199  & [4,4]    & 6  & 8  &       \\   \hline

   210  & [4,4]    & 6  & 8  &       \\   \hline

   211  & [4,2]    & 4  & 6  &       \\   \hline

   212  & [4,4]    & 6  & 8  &       \\   \hline

   222  & [4,4]    & 6  & 8  &       \\   \hline

   227  & [4,4]    & 6  & 8  &       \\   \hline

   233  & [4,4]    & 6  & 8  &       \\   \hline

   235  & [4,2]    & 4  & 6  &       \\   \hline

   236  & [4,2]    & 4  & 6  &       \\   \hline

   240  & [4,4]    & 6  & 8  &       \\   \hline
      & & & & \\[-.12in]
   267  & [4]    & 3  & 5  &  $Q_{64}$     \\   \hline
\end{tabular}
\end{center}

\newpage

\begin{center}
{\bf Table 3: Indecomposable, non 2--central,  2--groups of order $\leq 32$}

\bigskip

\begin{tabular}{|c|c|c|c|c|c|c|c|c|c|} \hline
 &&&& & & & & & \\[-.12in]
Order & \# & Type & Depth & Rank & $e(G)$ & $e^{\prime}(G)$ & $d_0(G)$ & $C_G(V)$'s & Notes \\   \hline
 &&&& & & & & & \\[-.12in]
8 & 4 &[2]& 2 & 2 & 1  &  -1   &  0 & ${\mathbf 2}^2$  &   $D_8$    \\   \hline
&&&&   & & & & & \\[-.12in]
16 & 8 &[4]& 1 & 2 & 3  &  -1   &  1 & ${\mathbf 4} \times {\mathbf 2}$  &   $AES_{16}$    \\   \hline
&&&&   & & & & & \\[-.12in]
 & 9 &[2,2]& 2 & 3 & 2  &  1   &  1 & ${\mathbf 2}^2$  &       \\   \hline
&&&&   & & & & & \\[-.12in]
 & 11 &[4]& 1 & 2 & 3 &  2 &  2 & ${\mathbf 4} \times {\mathbf 2}$  &       \\   \hline
&&&&   & & & & & \\[-.12in]
 & 12 &[2]& 1 & 2 & 1  &  -1   & 0 & ${\mathbf 2}^2$  &   $D_{16}$    \\   \hline
&&&&   & & & & & \\[-.12in]
 & 13 &[4]& 1 & 2 & 3 &  2 & 2 & ${\mathbf 2}^2$  &  $SD_{16}$     \\   \hline
&&&&   & & & & & \\[-.12in]
32 & 16 &[4,2]& 2 & 3 & 4 &  3 &  3 & ${\mathbf 4} \times {\mathbf 2}^2$  &       \\   \hline
&&&&   & & & & & \\[-.12in]
 & 17 &[4]& 2 & 2 & 3 &  -1 & 1 & ${\mathbf 8} \times {\mathbf 2}$  &       \\   \hline
&&&&   & & & & & \\[-.12in]
 & 20 &[2,2]& 2 & 3 & 2 &  1 &  1 & ${\mathbf 4} \times {\mathbf 2}^2$  &       \\   \hline
&&&&   & & & & & \\[-.12in]
 & 22 &[4]& 1 & 2 & 3 &  2 &  2 & ${\mathbf 8} \times {\mathbf 2}$  &       \\   \hline
&&&&   & & & & & \\[-.12in]
 & 26 &[4]& 2 & 2 & 3 &  -1 & 1 & ${\mathbf 8} \times {\mathbf 2}$, ${\mathbf 4} \times {\mathbf 2}$  &       \\   \hline
&&&&   & & & & & \\[-.12in]
 & 27 &[2,2]& 2 & 3 & 2 & 1 &  1 & ${\mathbf 2}^3$  &       \\   \hline
&&&&   & & & & & \\[-.12in]
 & 31 &[4]& 2 & 2 & 3 &  -1 & 2 & ${\mathbf 4} \times {\mathbf 4}$, ${\mathbf 4} \times {\mathbf 2}$  &       \\   \hline
&&&&   & & & & & \\[-.12in]
 & 32 &[4]& 1 & 2 & 3 &  2 &  2 & ${\mathbf 8} \times {\mathbf 2}$  &       \\   \hline
&&&&   & & & & & \\[-.12in]
 & 33 & [2,2] & 3 & 4 & 2 & -1 &  0 & ${\mathbf 2}^4$, ${\mathbf 2}^3$  &       \\   \hline
&&&&   & & & & & \\[-.12in]
 & 34 &[2,2]& 3 & 3 & 2 &  -1 & 0 & ${\mathbf 2}^3$  &       \\   \hline
&&&&   & & & & & \\[-.12in]
 & 36 &[2,2]& 3 & 3 & 2 &  -1 & 1 & ${\mathbf 4} \times {\mathbf 2}^2$, ${\mathbf 2}^3$  &       \\   \hline
&&&&   & & & & & \\[-.12in]
 & 37 &[4,2]& 2 & 3 & 4 &  3 &  3 & $\mathbf 4 \times {\mathbf 2}^2$  &       \\   \hline
&&&&   & & & & & \\[-.12in]
 & 38 &[4,2]& 2 & 3 & 4 &  2 &  2 & $\mathbf 4 \times {\mathbf 2}^2$, ${\mathbf 2}^3$  &       \\   \hline
&&&&   & & & & & \\[-.12in]
 & 39 &[4,2]& 2 & 3 & 4 &  3 &  3 & ${\mathbf 2}^3$  &       \\   \hline
&&&&   & & & & & \\[-.12in]
 & 41 &[4,4]& 2 & 3 & 6 &  5 &  5 & ${\mathbf 2}^3$  &       \\   \hline
&&&&   & & & & & \\[-.12in]
 & 42 &[4]& 3 & 3 & 3 &  -1 & 0 & ${\mathbf 2}^3$  &  $D_8*D_8$    \\   \hline
&&&&   & & & & & \\[-.12in]
 & 43 &[8]& 2 & 2 & 7 &  -1 & 3 & $Q_8 \times \mathbf 2$  &  $D_8*Q_8$     \\   \hline
&&&&   & & & & & \\[-.12in]
 & 44 &[4]& 2 & 3 & 3 & -1 &  1 & $\mathbf 4 \times \mathbf 2$, ${\mathbf 2}^3$  &       \\   \hline
&&&&   & & & & & \\[-.12in]
 & 45 &[8]& 1 & 2 & 7 &  4 &  4 & $Q_8 \times \mathbf 2$, $\mathbf 4 \times \mathbf 2$  &       \\   \hline
&&&&   & & & & & \\[-.12in]
 & 46 &[4]& 2 & 3 & 3 &  -1 &  3 & $Q_8 \times \mathbf 2$, ${\mathbf 2}^3$  &       \\   \hline
&&&&   & & & & & \\[-.12in]
 & 47 &[4]& 1 & 3 & 3 &  1 &  1 & $D_8 \times \mathbf 2$, ${\mathbf 2}^3$  &       \\   \hline
&&&&   & & & & & \\[-.12in]
 & 48 &[8]& 1 & 2 & 7 &  6 &  6 & $Q_8 \times \mathbf 2$  &       \\   \hline
&&&&   & & & & & \\[-.12in]
 & 49 &[2]& 2 & 2 & 1 &  -1 & 0 & ${\mathbf 2}^2$  & $D_{32}$      \\   \hline
&&&&   & & & & & \\[-.12in]
& 50 & [4] & 1& 2 & 3  & 2 & 2 & ${\mathbf 2}^2$  &  $SD_{32}$     \\   \hline
\end{tabular}
\end{center}

\newpage

\begin{center}
{\bf Table 4: Indecomposable, non 2--central,  order $64$, with $Cess^*(G) \neq \bold 0$}

\bigskip

\begin{tabular}{|c|c|c|c||c|c|c|c||c|c|c|c|} \hline
 &&&&&& & & & & & \\[-.12in]
\# & Type & Rank & $e^{\prime}(G)$ & \# & Type & Rank & $e^{\prime}(G)$ & \# & Type & Rank & $e^{\prime}(G)$\\   \hline
42 & [4] & 2& 2    &173 & [4,4] & 4&3     &193 & [4,2] & 3&3    \\   \hline
67 & [4] & 2& 2    &175 & [4,4] & 4&3     &196 & [4,2] & 3&2    \\   \hline
143 & [4] & 2&2     &183 & [4,4] & 4&5     &197 & [4,2] & 3&3    \\   \hline
182 & [4] & 2&2     &202 & [4,2] & 4&2     &198 & [4,2] & 3&3    \\   \hline
245 & [8] & 2&4     & &  & &              &200 & [4,4] & 3&5    \\   \hline
246 & [4] & 2&2     &32 & [4,2] & 3&3     &204 & [4,2] & 3&3    \\   \hline
249 & [8] & 2&6     &33 & [4,2] & 3&3     &206 & [4,2] & 3&3    \\   \hline
255 & [8] & 2&6     &40 & [2,2] & 3&1     &207 & [4,2] & 3&3    \\   \hline
258 & [8] & 2&4     &54 & [4,2] & 3&3     &208 & [4,2] & 3&3    \\   \hline
266 & [4] & 2&2     &60 & [4,2] & 3&3    &209 & [4,2] & 3&3    \\   \hline
 &  & &             &61 & [4,2] & 3&3     &213 & [4,2] & 3&2    \\   \hline
121 & [8] & 3&3     &62 & [2,2] & 3&1     &214 & [4,2] & 3&2    \\   \hline
130 & [8] & 3&3     &79 & [4,4] & 3&5     &215 & [4,4] & 3&5    \\   \hline
133 & [8] & 3&4     &80 & [4,4] & 3&4     &216 & [4,4] & 3&5    \\   \hline
180 & [8] & 3&3     &95 & [4,2] & 3&3     &218 & [4,2] & 3&3    \\   \hline
181 & [8] & 3&3     &97 & [4,2] & 3&3    &219 & [4,2] & 3&2    \\   \hline
247 & [4] & 3&1     &98 & [4,2] & 3&3     &220 & [4,2] & 3&3    \\   \hline
251 & [8] & 3&4     &99 & [4,2] & 3&2     &221 & [4,4] & 3&5    \\   \hline
253 & [8] & 3&3     &100 & [4,2] & 3&3     &223 & [4,4] & 3&5    \\   \hline
254 & [8] & 3&4     &102 & [4,4] & 3&5     &224 & [4,4] & 3&5    \\   \hline
257 & [8] & 3&3     &108 & [8,2] & 3&7     &225 & [4,2] & 3&2    \\   \hline
262 & [8] & 3&3     &115 & [8,2] & 3&7     &226 & [4,2] & 3&3    \\   \hline
 &  & &             &116 & [4,2] & 3&3     &228 & [4,2] & 3&2    \\   \hline
81 & [2,2,2] & 5&1     &118 & [4,2] & 3&3     &229 & [4,2] & 3&3    \\   \hline
 &  & &             &129 & [4,2] & 3&3     &230 & [4,2] & 3&3    \\   \hline
83 & [2,2,2] & 4&2     &132 & [4,2] & 3&3     &231 & [4,4] & 3&5    \\   \hline
85 & [2,2,2] & 4&2     &138 & [2,2] & 3&1     &232 & [4,4] & 3&5    \\   \hline
86 & [2,2,2] & 4&2     &161 & [4,4] & 3&4     &234 & [2,2] & 3&1    \\   \hline
89 & [4,2,2] & 4&4     &165 & [4,4] & 3&4     &238 & [4,4] & 3&5    \\   \hline
91 & [2,2,2] & 4&2     &166 & [4,4] & 3&4     &239 & [4,2] & 3&3    \\   \hline
146 & [2,2,2] & 4&2     &167 & [4,4] & 3&4     & &  & &    \\   \hline
147 & [4,2,2] & 4&4     &168 & [4,4] & 3&5     & &  & &    \\   \hline
148 & [2,2,2] & 4&2     &172 & [8,2] & 3&5     & &  & &    \\   \hline
150 & [2,2,2] & 4&2     &174 & [4,4] & 3&5     & &  & &    \\   \hline
151 & [2,2,2] & 4&2     &177 & [4,4] & 3&3     & &  & &    \\   \hline
 &  & &                 &178 & [4,4] & 3&4     & &  & &    \\   \hline
94 & [4,2] & 4&2     &179 & [4,4] & 3&5     & &  & &    \\   \hline
113 & [4,2] & 4&2     &185 & [4,4] & 3&3     & &  & &    \\   \hline
131 & [4,2] & 4&2     &186 & [4,4] & 3&4     & &  & &    \\   \hline
163 & [4,4] & 4&3     &189 & [4,2] & 3&3     & &  & &    \\   \hline
\end{tabular}
\end{center}

\end{document}